\documentstyle[12pt]{article}
\input amssym.def
\input amssym.tex
\parskip=\medskipamount
\newtheorem{definition}{Definition}[section]
\newtheorem{lemma}[definition]{Lemma}
\newtheorem{theorem}[definition]{Theorem}
\newtheorem{proposition}[definition]{Proposition}
\newtheorem{corollary}[definition]{Corollary}
\newtheorem{remark}[definition]{\it Remark}
\newtheorem{example}[definition]{\it Example}

 1   
\font\ddpp=msbm10  scaled \magstep 1  
\def\R{\hbox{\ddpp R}}               
\def\N{\hbox{\ddpp N}}               
\def\C{\hbox{\ddpp C}}    

\def\lcf{\lbrack\! \lbrack}
\def\rcf{\rbrack\! \rbrack}

\begin{document}
\baselineskip=.55cm
\title{Duality and modular class of a Nambu-Poisson structure}
\author{R. Ib\'a\~nez$^{1}$, M. de Le\'on$^2$, B. L\'opez$^3$ J.C.
Marrero$^{4},$ E. Padr\'on$^{4} $
\\[10pt]
{\small\it $^1$Departamento de Matem\'aticas, Facultad de
Ciencias}\\[-8pt] {\small \it Universidad del Pa{\'\i}s
Vasco}\\[-8pt] {\small \it Apartado 644, 48080 Bilbao,
SPAIN}\\[-8pt] {\small \it E-mail: mtpibtor@lg.ehu.es}\\[8pt]
{\small \it $^2$Laboratory of Dynamical Systems, Mechanics and
Control}\\[-8pt] {\small \it Instituto de Matem\'aticas y
F{\'\i}sica Fundamental}\\[-8pt] {\small\it Consejo Superior de
In\-ves\-ti\-ga\-cio\-nes Ci\-en\-t{\'\i}\-fi\-cas} \\[-8pt]
{\small\it Serrano 123, 28006 Madrid, SPAIN}\\[-8pt] {\small\it
E-mail: mdeleon@imaff.cfmac.csic.es} \\[8pt] {\small\it
$^3$Departamento de Matem\'aticas Edificio de Matem\'aticas e
Inform\'atica}\\[-8pt] {\small\it Campus Universitario de Tafira,
Universidad de Las Palmas de Gran Canaria} \\[-8pt] {\small\it
35017 Las Palmas, Canary Islands, SPAIN}\\[-8pt] {\small\it
E-mail: blopez@dma.ulpgc.es}\\[8pt] {\small\it $^4$Departamento de
Matem\'atica Fundamental, Facultad de Matem\'aticas}\\[-8pt]
{\small\it Universidad de la Laguna, La Laguna} \\[-8pt]
{\small\it Tenerife, Canary Islands, SPAIN}\\[-8pt] {\small\it
E-mail: jcmarrer@ull.es, mepadron@ull.es } }
\date{\empty}

\maketitle
\begin{abstract}
In this paper we introduce cohomology and homology theories for
Nambu-Poisson manifolds. Also we study the relation between  the
existence of a duality for  these theories and the vanishing of a
particular Nambu-Poisson cohomology class, the modular class. The
case of a regular Nambu-Poisson structure and some singular
examples are discussed.
\end{abstract}
\begin{quote}


\end{quote}

\newpage

\setcounter{section}{0}
\section{Introduction}
\setcounter{equation}{0}

Homology and cohomology theories have shown to be good tools in
the study of Poisson geometry, as they have been in other areas of
geometry and physics. In particular, a lot of work has been done in
the study of Poisson cohomology and Poisson homology (see for
example \cite{V3,W}). Poisson cohomology (also known as
Lichnerowicz-Poisson cohomology) of a Poisson manifold $M$ was
introduced by Lichnerowicz \cite{Lich} as the cohomology of the
subcomplex of the Chevalley-Eilenberg complex of the Lie algebra
$C^{\infty}(M,\R)$ consisting of the 1-differentiable cochains
that are derivations in each argument with respect to the usual
product of functions. Poisson cohomology provides a good framework
to express deformation and quantization obstructions. On the other
hand, Poisson homology (also known as canonical homology) was
defined as the homology of the operator boundary $\delta$ on
differential forms considered geometrically by Koszul \cite{Ko}
and algebraically by Brylinski \cite{Br} by taking the classical
limit of the Hochschild boundary operator for a quantized Poisson
algebra. The notion of Poisson (resp. symplectic) harmonicity has
appeared also to be very interesting. These cohomology and
homology theories can be extended to Lie algebroids, which are
algebraic structures of great interest in mathematics and physics
\cite{Ma}. Lie algebroids are a generalization of Lie algebras and
tangent bundles and each Poisson manifold has associated a Lie
algebroid in a natural way. Recently, a Poincar\'e type duality
between cohomology and homology theories has been proved by
Evens, Lu and Weinstein \cite{ELW} and Xu \cite{Xu} by using the
modular class of the Poisson structure \cite{W0}. This special
Poisson cohomology class has also been used for the classification
of quadratic Poisson structures \cite{LX}.

The aim of this paper is to introduce similar cohomology
and homology theories for Nambu-Poisson structures, as well as the
study of a Poincar\'e type duality. The concept of a Nambu-Poisson
structure was given by Takhtajan \cite{Ta} in 1994 in order
to find an axiomatic formalism for the $n$-bracket operation
\[
\{f_1,\dots ,f_n\}=\mbox{\rm det}(\frac{\partial f_i}{\partial x_j})
\]
proposed by Nambu \cite{N} and picking up the idea that in
statistical mechanics the basic result is Liouville theorem, which
follows from but does not require hamiltonian dynamics. A
Nambu-Poisson manifold of order $n$ is a manifold $M$ endowed with
a skew-symmetric $n$-bracket of functions $\{\,, \dots, \,\}$
satisfying the Leibniz rule and the fundamental identity
\[
\{f_1, \dots, f_{n-1}, \{g_1, \dots, g_n\} \} =
\sum_{i=1}^{n} \{g_1, \dots, \{f_1, \dots, f_{n-1}, g_i\}, \dots, g_n
\} ,
\]
for all $f_1, \dots, f_{n-1}, g_1, \dots, g_{n}$ $C^\infty$
real-valued functions on $M$. Note that the $n$-bracket
$\{\,,\dots ,\,\}$ allows us to introduce the Nambu-Poisson
$n$-vector $\Lambda$ characterized by the relation
$\Lambda(df_1,\dots ,df_n)=\{f_1,\dots ,f_n\}.$ The structure is
said to be regular if $\Lambda\not=0$ at every point. Recently,
local and global properties of Nambu-Poisson manifolds have been
studied \cite{AlGu,Gau,GM,ILMM,MVV,Nak,V4}. The canonical example
of a Nambu-Poisson structure  of order $n$ greater than 2 is the
one induced by a volume form on an oriented manifold of dimension
$n.$ In fact, a Nambu-Poisson manifold of order $n,$ $n\geq 3$,
admits a generalized foliation (the characteristic foliation)
whose leaves are either points or $n$-dimensional manifolds
endowed with a volume Nambu-Poisson structure. A strong effort is
being done in order to understand the geometry of Nambu-Poisson
structures, and also to know the Nambu mechanics (see for example
\cite{ChT,Cz,DFST}).

Recently, the authors have defined in \cite{ILMP} the notion of a
Leibniz algebroid in the same way as for the case of a Lie
algebroid but taking in mind the concept of Leibniz algebra
\cite{L1,L2}. A Leibniz algebra is a real vector space ${\frak g}$
endowed with a $\R$-bilinear mapping $\{\;,\;\}$ satisfying the
Leibniz identity
\[
\{a_1,\{a_2,a_3\}\}-\{\{a_1,a_2\}, a_3\}-\{a_2,\{a_1,a_3\}\}=0,
\]
for $a_1,a_2,a_3\in {\frak g}.$ If the bracket is skew-symmetric
we recover the notion of Lie algebra. In \cite{ILMP}, it was shown
that each Nambu-Poisson manifold $(M,\Lambda)$  of order $n$, with
$n\geq 3$, has associated a Leibniz algebroid, consisting in the
vector bundle $\Lambda^{n-1}(T^*M)\longrightarrow M$ whose space
of sections $ \Omega^{n-1}(M)$ has a Leibniz algebra structure
with bracket
\[
\lcf \alpha,\beta \rcf = {\cal L}_{\#_{n-1}(\alpha)}\beta+
(-1)^n(i(d\alpha)\Lambda)\beta,
\]
and a vector bundle homomorphism $\#_{n-1}: \Lambda^{n-1}(T^*M)
\longrightarrow TM$ given by $\#_{n-1}(\beta) = i(\beta)\Lambda$,
which provides a Leibniz algebra homomorphism between the spaces
of sections.
 The Leibniz algebroid
$(\Lambda^{n-1}(T^*M), \lcf\;,\;\rcf,\#_{n-1})$ allows us to
introduce the Leibniz algebroid cohomology. However, this
cohomology has infinite degrees and thus a Poincar\'e type
duality, with some homology theory, is not possible.

In this paper, in order to obtain a cohomology theory for
Nambu-Poisson manifolds without the problems above mentioned, we
begin by showing in Section $3$ a Lie algebra structure associated
with a Nambu-Poisson manifold $(M,\Lambda).$ In fact, we prove
that the center of the Leibniz algebra $(\Omega^{n-1}(M),\lcf
\;,\;\rcf)$ is the $C^\infty(M,\R)$-module
$\ker\#_{n-1}=\{\alpha\in \Omega^{n-1}(M)/\#_{n-1}(\alpha)=0\}$
and thus the quotient space
$\displaystyle\frac{\Omega^{n-1}(M)}{\ker\#_{n-1}}$ is a Lie
algebra. Moreover, if the Nambu-Poisson structure is regular,
$\displaystyle\frac{\Omega^{n-1}(M)}{\ker \#_{n-1}}$ is the space
of sections of the vector bundle
$\displaystyle\frac{\Lambda^{n-1}(T^*M)}{\ker\#_{n-1}}\rightarrow
M$ and this is a Lie algebroid. As a consequence of the above
results, we introduce in Section $4$ a cohomology theory for a
Nambu-Poisson manifold $(M,\Lambda)$. The resultant cohomology,
called Nambu-Poisson cohomology, is defined as the cohomology of
the Lie algebra
$\displaystyle\frac{\Omega^{n-1}(M)}{\ker\#_{n-1}}$ relative to a
certain representation. If the structure is regular, the
Nambu-Poisson cohomology is just the Lie algebroid cohomology of
$\displaystyle\frac{\Lambda^{n-1}(T^*M)}{\ker\#_{n-1}}\rightarrow
M$. So, we can think that for a Nambu-Poisson structure there
exists associated a kind of ``singular" Lie algebroid structure
and the corresponding cohomology. Also in Section $4$, we observe
that the characteristic foliation of a Nambu-Poisson manifold
allows us to introduce  the foliated cohomology which, in the
regular case, coincides with the  usual foliated cohomology
defined for regular foliations \cite{K,HMM,V1,V2}. Furthermore, in
this last case, we prove that  the foliated cohomology is
isomorphic to the Nambu-Poisson cohomology.

Section $5$ is devoted to the introduction of the canonical
Nambu-Poisson homology on an oriented Nambu-Poisson manifold. If
$M$ is an oriented manifold one can consider, in a natural way, a
homology complex whose $k$-chains are the $k$-vectors on $M$, the
homology operator on vector fields is the divergence with respect
to a volume and the resultant homology is dual of the de Rham
cohomology. The canonical Nambu-Poisson homology complex of an
oriented Nambu-Poisson manifold $(M,\Lambda)$ is a subcomplex of
this homology complex. In fact, if $(M,\Lambda)$ is regular, the
$k$-chains in the canonical Nambu-Poisson homology complex are the
$k$-vectors on $M$ which are tangent to the characteristic
foliation.

In Section $6$, we study the relation between the vanishing of the
modular class of an oriented Nambu-Poisson manifold $(M,\Lambda)$
and the existence of a duality between the homology and cohomology
theories introduced in the above sections. The modular class of
$M$ was defined in \cite{ILMP} by the authors and it was shown to
be null in some neighborhood of any regular point. An example of a
singular Nambu-Poisson structure with non null modular class was
also exhibited in \cite{ILMP}. Now, if $M$ is an oriented regular
Nambu-Poisson manifold of order $n$ $(n\geq 3)$ then, in Section
$6,$ we prove that the modular class of $M$ is null if and only if
there exists a basic volume with respect to the characteristic
foliation. Using this result, we obtain some interesting examples
of regular Nambu-Poisson structures with non null modular class.
Next, we show that the vanishing of the modular class implies the
existence of a duality between the foliated cohomology of $M$ and
the homology of a subcomplex of the canonical Nambu-Poisson
homology complex of $M.$ Thus, if $(M,\Lambda)$ is regular and
there exists a basic volume with respect to the characteristic
foliation of $M,$ we conclude that there is a duality between the
Nambu-Poisson cohomology and the canonical Nambu-Poisson homology
of $M.$

Finally, in Section $7$, we study a particular example, namely, a
singular Nambu-Poisson structure of order $3$ on $\R^3.$ We prove
that there is no duality between the canonical Nambu-Poisson
homology and the Nambu-Poisson cohomology  and that this last
cohomology is not isomorphic to the foliated cohomology.

\section{Preliminaries}
\setcounter{equation}{0}
All the manifolds considered in this paper are assumed to be connected.
\subsection{Nambu-Poisson structures}
Let $M$ be a differentiable manifold of dimension $m$. Denote by
${\frak X}(M)$ the Lie algebra of vector fields on $M$, by
$C^\infty(M,\R)$ the algebra of $C^\infty$ real-valued functions
on $M$, by $\Omega^{k}(M)$ the space of $k$-forms on $M$ and by
${\cal V}^k(M)$ the space of $k$-vectors.

{\it A Nambu-Poisson bracket} of order $n$ $(n\leq m)$ on $M$ (see
\cite{Ta}) is an $n$-linear mapping
$\{\;,\dots,\;\}:C^\infty(M,\R)\times \dots^{(n}\dots  \times
C^\infty(M,\R)\rightarrow C^\infty(M,\R)$ satisfying the following
properties :
\begin{enumerate}
\item[$(1)$] {\it Skew-symmetry:}
\[
\{f_1,\dots ,f_n\}=(-1)^{\varepsilon(\sigma)}\{f_{\sigma(1)},\dots
,f_{\sigma(n)}\},
\]
 for all $f_1,\dots ,f_n\in C^\infty(M,\R)$
and $\sigma\in Symm(n)$, where $Symm(n)$ is a symmetric group of
$n$ elements and $\varepsilon(\sigma)$ is the parity of the
permutation $\sigma$.
\item[$(2)$] {\it Leibniz rule:}
\[
\{f_1g_1, f_2,\dots ,f_n\}=f_1\{g_1,f_2,\dots
,f_n\}+g_1\{f_1,f_2,\dots ,f_n\},
\]
for all $f_1,\dots ,f_{n},g_1\in C^\infty(M,\R).$
\item[$(3)$] {\it Fundamental identity}:
\[
\{f_1,\dots ,f_{n-1},\{g_1,\dots ,g_n\}\}=
\sum_{i=1}^{n}\{g_1,\dots ,\{f_1,\dots ,f_{n-1},g_i\},\dots ,g_n\}
\]
for all $f_1,\dots ,f_{n-1},g_1,\dots ,g_{n}$ functions on $M.$
\end{enumerate}
Given a Nambu-Poisson bracket, we can define a skew-symmetric tensor
$\Lambda$ of type $(n,0)$
($n$-vector) as follows
\[
\Lambda(df_1,\dots ,df_n)=\{f_1,\dots,f_n\},
\]
for $f_1,\dots ,f_n\in C^\infty(M,\R).$ The pair $(M,\Lambda)$ is
called {\it  a Nambu-Poisson manifold of order $n.$ }

Let $(M,\Lambda)$ be a Nambu-Poisson manifold of order $n$ and $k$
be an integer with $k\leq n.$

If $\Lambda^k(T^*M)$ (respectively, $\Lambda^{n-k}(TM))$ denotes
the vector bundle of the $k$-forms (respectively, $(n-k)$-vectors)
then $\Lambda$ induces a homomorphism of vector bundles
$\#_k:\Lambda^k(T^*M)\rightarrow \Lambda^{n-k}(TM)$ by defining
\begin{equation}\label{0.2}
\#_k(\beta)=i(\beta)\Lambda(x),
\end{equation}
for $\beta\in \Lambda^k(T^*_xM)$ and $x\in M,$ where $i(\beta)$ is
the contraction by $\beta$. Denote also by $\#_k$ the homomorphism
of $C^\infty(M)$-modules from the space $\Omega^k(M)$ onto the
space  ${\cal V}^{n-k}(M)$ given by
\begin{equation}\label{0.3}
\#_k(\alpha)(x)=\#_k(\alpha(x)),
\end{equation}
for all $\alpha\in \Omega^{k}(M)$ and $x\in M.$
\begin{remark}{\rm It is clear that the mapping
$\#_k:\Omega^k(M)\rightarrow {\cal V}^{n-k}(M)$ induces an
isomorphism of $C^\infty(M,\R)$-modules
$\overline{\#_k}:\displaystyle\frac{\Omega^k(M)}{\ker
\#_k}\rightarrow \#_k(\Omega^k(M))$ defined by
\begin{equation}\label{2.2'}
\overline{\#_k}([\alpha])=\#_k(\alpha),
\end{equation}
for $[\alpha]\in \displaystyle\frac{\Omega^k(M)}{\ker\#_k}.$}
\end{remark}

If $f_1,\dots, f_{n-1}$ are $n-1$ functions on $M$, we define a
vector field
\begin{equation}\label{n2.3'}
X_{f_1\dots f_{n-1}}=\#_{n-1}(df_1\wedge \dots \wedge df_{n-1})
\end{equation}
which is called the {\it Hamiltonian vector field } associated with
the Hamiltonian functions $f_1,\dots ,f_{n-1}$.

>From the fundamental identity, it follows that the Hamiltonian
vector fields are infinitesimal automorphisms of $\Lambda,$ i.e.,
\begin{equation}\label{nn}
{\cal L}_{X_{f_1\dots f_{n-1}}}\Lambda=0,
\end{equation}
for all $f_1,\dots ,f_{n-1}\in C^\infty(M,\R).$

\begin{example}\label{volume}{\rm Let $M$ be an oriented $m$-dimensional manifold
and choose a volume form $\nu_M$ on $M$. Then, we can consider the
following Nambu-Poisson bracket $\{,\dots ,\}$ defined by the
formula
\[
df_1\wedge \dots \wedge df_{m}=\{f_1,\dots ,f_{m}\}{\nu_M}.
\]
In this case the homomorphisms $\#_k$ are isomorphisms, for all
$k\leq m$ (see
\cite{Gau}).}
\end{example}

The following theorem describes the local structure of the
Nambu-Poisson brackets of order $n$, with $n\geq 3$.

\begin{theorem}\label{2.2}\cite{AlGu,Gau,ILMM,MVV,Nak}
Let $M$ be a differentiable manifold of dimension $m$. The
$n$-vector $\Lambda$, $n \geq 3$, defines a Nambu-Poisson
structure on $M$ if and only if for all $x\in M$ with
$\Lambda(x)\not=0,$ there exist local coordinates $(x^1,\dots
,x^n,x^{n+1},\dots ,x^m)$ around $x$ such that
\[
\Lambda=\frac{\partial}{\partial x^1}\wedge \dots \wedge
\frac{\partial}{\partial x^n}.
\]
\end{theorem}

A point $x$ of a Nambu-Poisson manifold $(M,\Lambda)$ of order
$n\geq 3$ is said to be {\it regular } if $\Lambda(x)\not=0.$ If
every point of $M$ is regular then the Nambu-Poisson manifold
$(M,\Lambda)$ is said to be  {\it regular}.

Let $(M,\Lambda)$ be a Nambu-Poisson manifold of order $n$, with
$n\geq 3,$ and consider {\it the characteristic distribution } ${\cal D}$ on
$M$, given by
\begin{equation}\label{dc}
\begin{array}{rcl}
x\in M\rightarrow {\cal D}(x)
&=&\#_{n-1}(\Lambda^{n-1}(T_x^*M))\\&=&<\{X_{f_1\dots
f_{n-1}}(x)/f_1,\dots ,f_{n-1}\in C^\infty(M,\R)\}>\subseteq T_xM.
\end{array}
\end{equation}
Then,  ${\cal D}$ defines a generalized foliation on $M$ whose
leaves are either points or $n$-di\-men\-sio\-nal manifolds
endowed with a Nambu-Poisson structure coming from a volume form
(see \cite{ILMM}).

\begin{remark}\label{r2.2'}{\rm Let $(M,\Lambda)$ be an $m$-dimensional regular
Nambu-Poisson manifold of order $n,$ with $n\geq 3.$ From Theorem
\ref{2.2}, we deduce:

\medskip

$(i)$ ${\cal D}$ defines a foliation on $M$ of dimension $n$.

\medskip

$(ii)$ For all $k\in \{0,\dots ,n\},$ $\ker\#_k$ (respectively,
$\#_k(\Lambda^k(T^*M))$) is a vector subbundle of
$\Lambda^k(T^*M)\rightarrow M$ (respectively,
$\Lambda^{n-k}(TM)\rightarrow M$) of rank {\small
$\left(\begin{array}{l}m\\k\end{array}\right)-
\left(\begin{array}{l}n\\k\end{array}\right)$} (respectively,
{\small $\left(\begin{array}{l}n\\k\end{array}\right)$}) and the
homomorphism $\#_k:\Lambda^k(T^*M)\rightarrow \Lambda^{n-k}(TM)$
induces an isomorphism of vector bundles
\[
\overline{\#_k}:\frac{\Lambda^k(T^*M)}{\ker \#_k}\rightarrow
\#_k(\Lambda^k(T^*M)).
\]
The notation $\overline{\#_k}$ is justified by the following fact.
The space of the $C^\infty$-differentiable sections of
$\displaystyle\frac{\Lambda^k(T^*M)}{\ker \#_k}\rightarrow M$
(respectively, $\#_k(\Lambda^k(T^*M))\rightarrow M)$ can be
identified with $\displaystyle\frac{\Omega^k(M)}{\ker \#_k}$
(respectively, $\#_k(\Omega^k(M))$) in such a sense that the
corresponding isomorphism of $C^\infty(M,\R)$-modules induced by
$\overline{\#_k}$ is just the mapping
$\overline{\#_k}:\displaystyle\frac{\Omega^k(M)}{\ker
\#_k}\rightarrow \#_k(\Omega^k(M))$ given by (\ref{2.2'}).

\medskip

$(iii)$ The $C^\infty$-differentiable sections of the vector
bundle $\#_k(\Lambda^k(T^*M))\rightarrow M$ are the
$(n-k)$-vectors on $M$ which are tangent to ${\cal D}.$ We recall
that an $(n-k)$-vector $P$ on $M$ is tangent to ${\cal D}$ if
\[
i(\alpha(x))(P(x))=0,
\]
for all $x\in M$ and for all $\alpha(x)\in {\cal D}^0(x),$ where
${\cal D}^0(x)$ is the annihilator of ${\cal D}(x)$ in $T_x^*M.$
Note that ${\cal D}^0(x)=\ker (\#_{1|T_x^*M}),$ for all $x\in M.$}
\end{remark}

\subsection{ The Leibniz algebroid associated with a
Nam\-bu\--Poi\-sson  struc\-tu\-re} In \cite{ILMP} we have
introduced the notion of a Leibniz algebroid, a natural
generalization of the notion of a Lie algebroid, and we have
proved that every Nambu-Poisson manifold has associated a
canonical Leibniz algebroid. Next, we will describe this
structure.

First, we recall the definition of real Leibniz algebra (see
\cite{Cu,L1,L2,LP}). A {\it Leibniz algebra structure}  on a real
vector space ${\frak g}$ is a $\R$-bilinear map $\{\;,\;\}:{\frak
g}\times {\frak g}\rightarrow {\frak g}$ satisfying the {\it
Leibniz identity}, that is,
\[
\{a_1,\{a_2,a_3\}\}-\{\{a_1,a_2\},a_3\}-\{a_2,\{a_1,a_3\}\}=0,
\]
for $a_1, a_2, a_3\in {\frak g}.$ In such a case, the pair
$({\frak g},\{\;,\;\})$ is called a {\it Leibniz algebra}.

Moreover, if the skew-symmetric
condition is required then $({\frak g}, \{\;,\;\})$ is a Lie algebra.
In this sense, a Leibniz algebra is a non-commutative version of a
Lie algebra.

The notion of Leibniz algebroid can be introduced in the same way as
that of Lie algebroid.

\begin{definition}\label{alg}
A Leibniz algebroid structure on a differentiable vector bundle
$\pi:E\rightarrow M$ is a pair that consists of a Leibniz algebra
structure $\lcf \;,\; \rcf $ on the space $\Gamma(E)$ of the global cross
sections of $\pi:E\longrightarrow M$ and a vector bundle morphism
$\varrho:E\rightarrow TM,$ called the anchor map, such that the
induced map $\varrho:\Gamma(E)\longrightarrow \Gamma(TM)={\frak
X}(M)$ satisfies the following relations:
\begin{enumerate}
\item
$\varrho\lcf s_1,s_2 \rcf =[\varrho(s_1),\varrho(s_2)],$
\item
$\lcf s_1,fs_2\rcf =f\lcf s_1,s_2\rcf +\varrho(s_1)(f)s_2,$
\end{enumerate}
for all $s_1,s_2\in \Gamma(E)$ and $f\in C^\infty(M,\R).$

A triple $(E,\lcf \;,\; \rcf,\varrho)$ is called a Leibniz algebroid
over $M$.
\end{definition}

Every Lie algebroid over a manifold $M$ is trivially a Leibniz
algebroid. In fact, a Leibniz algebroid
$(E,\lcf\;,\;\rcf,\varrho)$ over $M$ is a Lie algebroid if and
only if the Leibniz bracket $\lcf\;,\;\rcf$ on $\Gamma(E)$ is
skew-symmetric.

Now, let $(M,\Lambda)$ be a Nambu-Poisson manifold of order $n,$
$n\geq 3,$ and
 ${\cal L}$ the Lie derivative operator on $M.$
The Leibniz algebroid attached  to $M$ is just the triple
$(\Lambda^{n-1}(T^*M),\lcf \;,\; \rcf,$ $\#_{n-1})$, where $\lcf
\;,\; \rcf:\Omega^{n-1}(M)\times
\Omega^{n-1}(M)\rightarrow\Omega^{n-1}(M)$ is the bracket of
$(n-1)$-forms defined by
\begin{equation}\label{0.1}
\lcf \alpha,\beta \rcf={\cal
L}_{\#_{n-1}(\alpha)}\beta+(-1)^n\#_{n}(d\alpha)\beta,
\end{equation}
for all $\alpha,\beta\in \Omega^{n-1}(M).$ In particular we have
that
\begin{equation}\label{sc}
\#_{n-1}(\lcf\alpha,\beta\rcf)=[\#_{n-1}(\alpha),\#_{n-1}(\beta)],
\end{equation}
for all $\alpha,\beta\in \Omega^{n-1}(M).$

Moreover, in \cite{ILMP} it was proved that the only non-null
Nambu-Poisson structures of order greater than two on an oriented
manifold $M$ of dimension $m$ such that its Leibniz algebroid is a
Lie algebroid are those defined by non-null $m$-vectors.

Let $(E,\lcf\;,\;\rcf,\varrho)$ be a Leibniz algebroid over a
manifold $M.$ For every $k\in \N$, we con\-si\-der the vector
space
\[
C^k(\Gamma(E);C^\infty(M,\R))=\{c^k:\Gamma(E)\times
\dots^{(k}\dots  \times \Gamma(E)\rightarrow C^\infty(M,\R)/ c^k
\mbox{ is $k$-linear}\}
\] and the operator
$\partial:C^k(\Gamma(E);C^\infty(M,\R))\rightarrow
C^{k+1}(\Gamma(E);C^\infty(M,\R))$ defined by
\[
\begin{array}{rcl}
\partial c^k(s_0,\dots
,s_k)&=&\displaystyle\sum_{i=0}^k(-1)^i\varrho(s_i)(c^k(s_0,\dots
,\widehat{s_i},\dots, s_k))\\[5pt] &&\kern-20pt
+\displaystyle\sum_{0\leq i<j\leq k}(-1)^{i-1}c^k(s_0,\dots
,\widehat{s_i},\dots ,s_{j-1},\lcf s_i,s_j\rcf, s_{j+1},\dots
,s_k),
\end{array}
\]
for $c^k\in C^k(\Gamma(E);C^\infty(M,\R))$ and $s_0,\dots ,s_k\in
\Gamma(E).$

Then, it follows that $\partial^2=0.$ The resultant cohomology is
called the {\it Leibniz algebroid cohomology of $E.$ } This
cohomology also can be described as the one defined by the
representation
\[
\Gamma(E)\times C^\infty(M,\R)\rightarrow C^\infty(M,\R),
\makebox[1cm]{} (s,f)\mapsto \varrho(s)(f).
\]
The definition of the cohomology of a Leibniz algebra relative to
a representation can found in \cite{L1,L2,LP}.

Note that if $c^k\in C^k(\Gamma(E);C^\infty(M,\R))$ is
skew-symmetric (respectively, $C^\infty(M,\R)$-linear) then, in
general, $\partial c^k$ is not skew-symmetric (respectively,
$C^\infty(M,\R)$-linear) (for more details, see \cite{ILMP}).

Nevertheless, if $(E,\lcf\;,\;\rcf, \varrho)$ is a Lie algebroid
and $c^k\in C^k(\Gamma(E);C^\infty(M,\R))$ is skew-symmetric and
$C^\infty(M,\R)$-linear then $\partial c^k$ is also skew-symmetric
and $C^\infty(M,\R)$-linear. Thus, in this case, we can consider
the subcomplex of $(C^*(\Gamma(E);$ $ C^\infty(M,\R)),\partial^*)$
that consists of the skew-symmetric $C^\infty(M,\R)$-linear
cochains. The cohomology of this subcomplex is just {\it the Lie
algebroid cohomology of $E$} (see \cite{Ma}).

\begin{remark}{\rm Let $(M,\Lambda)$ be a Nambu-Poisson manifold
of order $n,$ with $n\geq 3,$ and $(\Lambda^{n-1}(T^*M),\lcf
\;,\;\rcf, \#_{n-1})$ the corresponding Leibniz algebroid. Now,
the Leibniz algebroid cohomology operator is given by
\begin{equation}\label{opal}
\begin{array}{rcl}
\partial c^k(\alpha_0,\dots
,\alpha_k)&=&\displaystyle\sum_{i=0}^k(-1)^i\#_{n-1}(\alpha_i)(c^k(\alpha_0,\dots
,\widehat{\alpha_i},\dots, \alpha_k))\\ &&\kern-20pt
+\displaystyle\sum_{0\leq i<j\leq k}(-1)^{i-1}c^k(\alpha_0,\dots
,\widehat{\alpha_i},\dots ,\alpha_{j-1},\lcf
\alpha_i,\alpha_j\rcf, \alpha_{j+1},\dots ,\alpha_k),
\end{array}
\end{equation}
for all $c^k\in C^k(\Omega^{n-1}(M);C^\infty(M,\R))$ and
$\alpha_0,\dots ,\alpha_k\in \Omega^{n-1}(M).$}
\end{remark}

\section{A Lie algebra associated with a Nambu-Poisson manifold}
\setcounter{equation}{0}

If $({\frak g},[\;,\;])$ is  a Leibniz algebra, we define its
{\it center,} $Z({\frak g}),$ as the kernel of the adjoint
representation
\[
ad:{\frak g}\rightarrow \mbox{ End }({\frak
g}),\makebox[1cm]{}x\mapsto [x,\;\cdot\;].
\]
It is easy to prove that ${\frak g}/Z({\frak g})$ endowed with the
induced bracket is a Lie algebra (see \cite{Cu}).

In the particular case of a Nambu-Poisson manifold $(M,\Lambda)$
of order $n\geq 3$, we have that the center of the Leibniz algebra
$(\Omega^{n-1}(M),\lcf \;,\; \rcf)$  is the space
\[
Z(\Omega^{n-1}(M))=\{\alpha\in \Omega^{n-1}(M)\;/\;\lcf
\alpha,\beta \rcf=0, \;\forall\beta\in \Omega^{n-1}(M)\}
\]
and that $(\Omega^{n-1}(M)/Z(\Omega^{n-1}(M)),\lcf \;,\;
\rcf^{\tilde{\;}})$ is a Lie algebra, where
\[
\lcf\; ,\;\rcf^{\tilde{\;}}:\Omega^{n-1}(M)/Z(\Omega^{n-1}(M))\times
\Omega^{n-1}(M)/Z(\Omega^{n-1}(M)) \rightarrow
\Omega^{n-1}(M)/Z(\Omega^{n-1}(M))
\]
is the bracket given by
\begin{equation}\label{3.0'}
\lcf[\alpha] ,[\beta]\rcf^{\tilde{\;}}=[\lcf\alpha ,\beta\rcf],
\end{equation}
for all $[\alpha], [\beta]\in \Omega^{n-1}(M)/Z(\Omega^{n-1}(M)).$

The next result gives an explicit description of the center of
$(\Omega^{n-1}(M), \lcf\; ,\;\rcf).$

\begin{proposition}
Let $(M,\Lambda)$ be an $m$-dimensional Nambu-Poisson manifold of
order $n,$ with $n\geq 3$. Then, the  center of the algebra
$(\Omega^{n-1}(M),\lcf\; ,\;\rcf)$ is the $C^\infty(M,\R)$-module
\[
\ker \#_{n-1}=\{\alpha\in \Omega^{n-1}(M)/\#_{n-1}(\alpha)=0\}.
\]
\end{proposition}
{\it Proof:} If $\alpha$ is an $(n-1)$-form on $M$ such that
$\#_{n-1}(\alpha)=0$ then, from  (\ref{0.1}), it follows that
\begin{equation}\label{3.0}
\lcf\alpha ,\beta\rcf=(-1)^n\#_{n}(d\alpha)\beta,
\end{equation}
for all  $\beta\in \Omega^{n-1}(M).$

On the other hand, using a result proved in \cite{ILMP} (see
relation $(3.3)$ in \cite{ILMP}), we have that
\[
0={\cal L}_{\#_{n-1}(\alpha)}\Lambda
=(-1)^n\#_{n}(d\alpha)\Lambda.
\]
Thus, we deduce that $\#_n(d\alpha)=0.$ Consequently, $\lcf\alpha
,\beta\rcf=0$ (see (\ref{3.0})).

Conversely, suppose that $\alpha$ is an $(n-1)$-form on $M$ such that
\begin{equation}\label{3.4}
\lcf\alpha,\beta\rcf=0, \makebox[2cm]{ for all} \beta\in \Omega^{n-1}(M).
\end{equation}
In order to prove that $\#_{n-1}(\alpha)(x)=0,$ for all $x\in M$, we
distinguish two cases:

\medskip

$(i)$ If $\Lambda(x)=0,$ it is obvious that $\#_{n-1}(\alpha)(x)=0.$

\medskip

$(ii)$ If $\Lambda(x)\not=0$ then, using Theorem \ref{2.2}, we
have that there exist local coordinates $(x^1,\dots
,x^n,x^{n+1},\dots ,x^m)$ in a connected open neighborhood $U$ of
$x$ such that
\begin{equation}\label{3.5'}
\Lambda=\frac{\partial}{\partial x^1}\wedge \dots \wedge
\frac{\partial }{\partial x^n}.
\end{equation}
Now, the $(n-1)$-form $\alpha$ on $U$ can be written as follows
\begin{equation}\label{36'}
\alpha=\sum_{i=1}^n(-1)^{n-i}\alpha_idx^1\wedge \dots \wedge
\widehat{dx^i}\wedge \dots \wedge dx^n + \alpha',
\end{equation}
where $\alpha_i\in C^\infty(U,\R)$ and $\alpha'$ is an
$(n-1)$-form on $U$ satisfying the condition
$\#_{n-1}(\alpha')=0.$

Note that on $U$
\begin{equation}\label{3.5}
\#_{n-1}(\alpha)=\sum_{i=1}^n\alpha_i\frac{\partial }{\partial
x^i}.
\end{equation}
On the other hand, from (\ref{sc}), (\ref{3.4}), (\ref{3.5'}) and
(\ref{3.5}), we obtain that, for all $j\in \{1,\dots ,n\}$,
\[
0=\#_{n-1}(\lcf\alpha,(-1)^{n-j}dx^1\wedge\dots \wedge
\widehat{dx^j}\wedge \dots \wedge
dx^n\rcf)=[\#_{n-1}(\alpha),\frac{\partial}{\partial
x^j}]=-\sum_{i=1}^n\frac{\partial\alpha_i}{\partial
x^j}\frac{\partial}{\partial x^i}.
\]
Consequently,
\begin{equation}\label{3.7}
\frac{\partial \alpha_i}{\partial x^j}=0,\makebox[2cm]{for all
}i,j\in \{1,\dots ,n\}.
\end{equation}
This implies that (see (\ref{3.5'}) and (\ref{36'})) on $U$, we
have
\begin{equation}\label{3.6}
\#_n(d\alpha)=0.
\end{equation}
Moreover, we shall see that $d\alpha_i=0,$ for all $i\in \{1,\dots
,n\}$. Indeed, consider the $(n-1)$-forms
$\beta=dx^1\wedge\dots\wedge \widehat{dx^j}\wedge \dots  \wedge
dx^n,$ for all $j$. Since $\lcf\alpha,\beta\rcf=0,$ using
(\ref{3.5}) and (\ref{3.6}), we obtain
\[
0=\lcf\alpha,\beta\rcf= {\cal
L}_{\#_{n-1}(\alpha)}\beta=\sum_{i=1}^nd\alpha_i\wedge
i(\frac{\partial }{\partial x^i})(dx^1\wedge\dots \wedge
\widehat{dx^j}\wedge \dots \wedge dx^n).
\]
Thus, $\displaystyle\frac{\partial \alpha_i}{\partial x^k}=0$ for
all $k\in \{n+1,\dots ,m\}$ and for all $i\in \{1,\dots ,n\}.$
 This fact and (\ref{3.7}) imply that $d\alpha_i=0,$
that is, $\alpha_i$ is a real constant, for all $i\in \{1,\dots
,n\}.$

Next, we will prove that $\alpha_i=0$ for all $i\in
\{1,\dots,n\}$. We consider the $(n-1)$-form $\beta'$ on $U$ given
by
\[
\beta' =x^jdx^1\wedge\dots \wedge \widehat{dx^j}\wedge \dots
\wedge dx^n.
\]

Using (\ref{0.1}), (\ref{3.5}), (\ref{3.6}) and the fact that
$\alpha_i$ is constant, we have that
\[0=\lcf\alpha,\beta'\rcf=\alpha_j dx^1\wedge
\dots \wedge \widehat{dx^j}\wedge \dots \wedge dx^n.
\]
Therefore,
\[
\alpha_j=0, \makebox[.5cm]{}\mbox{for all }j\in \{1,\dots ,n\}.
\]
Finally, from  (\ref{3.5}) we conclude that $\#_{n-1}(\alpha)=0$
on $U.$ In particular,
\[
\#_{n-1}(\alpha)(x)=0.
\]
\hfill$\Box$

Hence, if $(M,\Lambda)$ is an $m$-dimensional Nambu-Poisson
manifold of order $n$, the quotient space
\[
\Omega^{n-1}(M)/Z(\Omega^{n-1}(M))=\Omega^{n-1}(M)/\ker \#_{n-1}
\]
is a $C^\infty(M,\R)$-module endowed with a skew-symmetric bracket
$\lcf\;,\;\rcf^{\tilde{\;}}$ given by (\ref{3.0'}) which satisfies
the Jacobi identity and the following property:
\begin{equation}
\lcf[\alpha],f[\beta]\rcf^{\tilde{\;}}=f\lcf[\alpha],[\beta]\rcf^{\tilde{\;}}
+\#_{n-1}(\alpha)(f)[\beta]
\end{equation}
for all $[\alpha],[\beta]\in \Omega^{n-1}(M)/\ker\#_{n-1}$ and
$f\in C^\infty(M,\R).$

Furthermore, using (\ref{sc}) we obtain that the mapping
$\widetilde{\#_{n-1}}:\Omega^{n-1}(M)/\ker\#_{n-1}\rightarrow
{\frak X}(M)$ defined by
\begin{equation}
\widetilde{\#_{n-1}}([\alpha])=\#_{n-1}(\alpha)
\end{equation}
induces a  homomorphism of Lie algebras between
$(\Omega^{n-1}(M)/\ker\#_{n-1},\lcf\;,\;\rcf^{\tilde{\;}})$ and
$({\frak X}(M),$ $[\;,\;]).$

\begin{remark}\label{r3.2}{\rm Let $(M,\Lambda)$ be a regular Nambu-Poisson
manifold of order $n,$ with $n\geq 3.$

\medskip

$(i)$ Using the above facts and Remark \ref{r2.2'}, we deduce that
the triple
\[
(\frac{\Lambda^{n-1}(T^*M)}{\ker\#_{n-1}},\lcf\;,\;\rcf^{\tilde{\;}},\widetilde{\#_{n-1}})
\]
is a Lie algebroid over $M.$

\medskip

$(ii)$ If ${\cal F}$ is a foliation on a manifold $N$ and
$F=\bigcup_{x\in N}{\cal F}(x)\rightarrow N$ is the corresponding
vector subbundle of $TN$ then the triple $(F,[\;,\;],i)$ is a Lie
algebroid over $N,$ where $[\;,\;]$ is the usual Lie bracket of
vector fields and $i:F\rightarrow TN$ is the inclusion.

\medskip

$(iii)$ If ${\cal D}$ is the characteristic foliation of $M,$ then
the Lie algebroids $(\bigcup_{x\in M}{\cal
D}(x)=\#_{n-1}(\Lambda^{n-1}(T^*M)),[\;,\;], i)$,
$(\displaystyle\frac{\Lambda^{n-1}(T^*M)}{\ker \#_{n-1}},
\lcf\;,\;\rcf^{\tilde{\;}}, \widetilde{\#_{n-1}})$ are isomorphic
(see Remark \ref{r2.2'}).}
\end{remark}

\section{The Nambu-Poisson cohomology and the foliated cohomology}
\setcounter{equation}{0} Let $(M,\Lambda)$ be a Nambu-Poisson
manifold of order $n, n\geq 3.$  According to the precedent
section, the quotient space
$\displaystyle\frac{\Omega^{n-1}(M)}{\ker \#_{n-1}}$ endowed with
the bracket $\lcf\;,\;\rcf^{\tilde{\;}}$ given by (\ref{3.0'}) is
a Lie algebra.

Moreover, using (\ref{sc}), we deduce that $C^\infty(M,\R)$ is a
$(\Omega^{n-1}(M)/\ker \#_{n-1})$-module relative to the
representation:
\[
\Omega^{n-1}(M)/\ker \#_{n-1}\times C^\infty(M,\R)\rightarrow
C^\infty(M,\R),\makebox[1cm]{} ([\alpha],f)\mapsto[\alpha]
f=(\#_{n-1}(\alpha))(f).
\]
Thus, one can consider the skew-symmetric cohomology complex
\[
\left(C^*(\Omega^{n-1}(M)/\ker\#_{n-1};C^\infty(M,\R))=\bigoplus_{k}C^k(\Omega^{n-1}(M)/\ker\#_{n-1};
C^\infty(M,\R)) ,\tilde\partial\right),
\]
where the space of the $k$-cochains $C^k(\Omega^{n-1}(M)/\ker
\#_{n-1};C^\infty(M,\R))$ consists of skew-symmetric
$C^\infty(M,\R)$-linear mappings
\[
c^k:(\Omega^{n-1}(M)/\ker\#_{n-1})\times \dots ^{(k}\dots \times
(\Omega^{n-1}(M)/\ker\#_{n-1})\rightarrow C^\infty(M,\R)
\]
and the cohomology operator $\tilde\partial$ is given by
\begin{equation}\label{partial}
\begin{array}{rcl}
\tilde\partial c^k([\alpha_0],\dots
,[\alpha_k])&=&\displaystyle\sum_{i=0}^k(-1)^i
(\#_{n-1}(\alpha_i))(c^k([\alpha_0],\dots
,\widehat{[\alpha_i]},\dots ,[\alpha_k]))\\ &&
\kern-70pt+\displaystyle\sum_{0\leq i<j\leq k}(-1)^{i-1}c^k(
[\alpha_0],\dots ,\widehat{[\alpha_i]},\dots
,[\alpha_{j-1}],[\lcf\alpha_i,\alpha_j\rcf], [\alpha_{j+1}],\dots
,[\alpha_k]),
\end{array}
\end{equation}
for all $c^k\in C^k(\Omega^{n-1}(M)/\ker\#_{n-1};C^\infty(M,\R)),$
and $[\alpha_0],\dots ,[\alpha_k]\in
\displaystyle\frac{\Omega^{n-1}(M)}{\ker\#_{n-1}}.$

The cohomology of this complex is called {\it the Nambu-Poisson
cohomology} and denoted by $H_{NP}^*(M).$

\begin{remark}\label{1-2}{\rm
Let $(M,\Lambda)$ be a Nambu-Poisson manifold of order $n,$ $n\geq
3.$ Consider \linebreak $(C^*(\Omega^{n-1}(M);$
$C^\infty(M,\R)),\partial)$ the cohomology complex associated with
the Leibniz algebroid $(\Lambda^{n-1}(T^*M),$
$\lcf\;,\;\rcf,\#_{n-1}).$ The natural projection
$p:\Omega^{n-1}(M)\rightarrow
\displaystyle\frac{\Omega^{n-1}(M)}{\ker\#_{n-1}}$ allows us to
define the  homomorphisms of $C^\infty(M,\R)$-modules
\[p^k:C^k(\Omega^{n-1}(M)/\ker\#_{n-1};C^\infty(M,\R))\rightarrow
C^k(\Omega^{n-1}(M);C^\infty(M,\R)),\makebox[1cm]{}c^k\mapsto
p^k(c^k),
\]
$p^k(c^k):\Omega^{n-1}(M)\times \dots^{(k}\dots \times
\Omega^{n-1}(M)\rightarrow C^\infty(M,\R)$ being  the mapping
given by
\[
p^k(c^k)(\alpha_1,\dots, \alpha_k)=c^k([\alpha_1],\dots
,[\alpha_k]).
\]
A direct computation, using (\ref{opal}) and (\ref{partial}),
proves that these homomorphisms induce a homomorphism between the
complexes $(C^*(\Omega^{n-1}(M)/\ker\#_{n-1};$
$C^\infty(M;\R)),\tilde\partial)$ and
$(C^*(\Omega^{n-1}(M);C^\infty(M,\R)),\partial)$. Therefore, we
have the corresponding homomorphism in cohomology
\[
p^*:H^*_{NP}(M)\rightarrow H^*(\Omega^{n-1}(M);C^\infty(M,\R)).
\]
Moreover, since the space of $0$-cochains in both complexes  is
$C^\infty(M,\R)$, then
\[
p^1:H^1_{NP}(M)\rightarrow H^1(\Omega^{n-1}(M);C^\infty(M,\R))
\]
is a monomorphism.}
\end{remark}

Now, using the isomorphism of $C^\infty(M,\R)$-modules
\begin{equation}\label{4.0'}
\overline{\#_{n-1}}: \Omega^{n-1}(M)/\ker \#_{n-1}\rightarrow
\#_{n-1}(\Omega^{n-1}(M)),\makebox[1cm]{}
\overline{\#_{n-1}}([\alpha])=\#_{n-1}(\alpha),
\end{equation}
we will relate the Nambu-Poisson cohomology with the foliated
cohomology of $(M, {\cal D})$, where ${\cal D}$ is the
characteristic foliation of $M$.

The foliated cohomology of $(M,{\cal D})$ is defined as follows.
We consider the space $\Omega^k(M,{\cal D})$  of the $k$-forms
$\alpha$ on $M$ such that
\[
\alpha(X_1,\dots ,X_k)=0,\mbox{ for all }X_1,\dots ,X_k\in
\#_{n-1}(\Omega^{n-1}(M)).
\]
>From (\ref{sc}), it follows that if $\alpha\in \Omega^k(M,{\cal
D})$ then $d\alpha\in \Omega^{k+1}(M,{\cal D}).$ Now, denote by
$\Omega^k({\cal D})$ the $C^\infty(M,\R)$-module
$\displaystyle\frac{\Omega^k(M)}{\Omega^k(M,{\cal D})}.$ Then, the
exterior differential induces a cohomology operator $
\tilde{d}:\Omega^k({\cal D})\rightarrow\Omega^{k+1}({\cal D})$
\begin{equation}\label{n4.2'}
\tilde{d}([\alpha])=[d\alpha],\makebox[.3cm]{} \mbox{for }[\alpha]
\in \Omega^k({\cal D}).
\end{equation}

The resultant cohomology $H^*({\cal D})$ is called the {\it foliated
cohomology} of $(M,{\cal D})$ and the operator $\tilde{d}$ is called the {\it foliated
differential} of $(M,{\cal D})$. Note that if $M$ is a regular
Nambu-Poisson manifold, $H^*({\cal D})$ is just the usual foliated
cohomology  of $(M,{\cal D})$ (see \cite{K,HMM,V1,V2}).

On the other hand, we have
\begin{proposition}\label{p4.0} Let $(M,\Lambda)$ be a
Nambu-Poisson manifold of order $n$, with $n\geq 3.$ Then,
\[
\Omega^k(M,{\cal D})=\ker\#_k,
\]
for all $k\in \{0,\dots ,n\}.$ Thus,
\[
\#_{k+1}(d\alpha)=0,
\]
for all $\alpha\in \ker \#_k.$
\end{proposition}
{\it Proof:} Suppose that $\alpha\in \Omega^k(M,{\cal D})$. We
will prove that $\#_{k}(\alpha)(x)=0,$ for all $x\in M.$

We distinguish two cases:

\medskip

$(i)$ If $\Lambda(x)=0,$ it is clear that $\#_k(\alpha)(x)=0.$

\medskip

$(ii)$ If $\Lambda(x)\not=0$ then, using Theorem \ref{2.2}, we
deduce that there exist local coordinates $(x^1,\dots
,x^n,x^{n+1},\dots ,x^m)$ in an open neighborhood $U$ of $x$ such
that
\[
\Lambda=\frac{\partial }{\partial x^1}\wedge \dots \wedge
\frac{\partial }{\partial x^n}.
\]
Now, we consider an $(n-1)$-form $\beta_i$ on $M$ satisfying
\[
\#_{n-1}(\beta_i)(x)=\frac{\partial }{\partial x^i}_{|x},
\]
for all $i\in \{1,\dots ,n\}.$ Since $\alpha\in \Omega^k(M,{\cal
D})$, it follows that
\[
\alpha(\#_{n-1}(\beta_{i_1}),\dots ,\#_{n-1}(\beta_{i_k}))=0,
\]
for all $1\leq i_1<\dots <i_{k}\leq n.$ Thus,
\[
\alpha_x(\frac{\partial }{\partial x^{i_1}}_{|x},\dots
,\frac{\partial }{\partial x^{i_k}}_{|x})=0.
\]
This implies that $\#_k(\alpha)(x)=0.$
Therefore, $\Omega^k(M,{\cal D})\subseteq \ker\#_k.$

The proof of the inclusion $\ker\#_k\subseteq \Omega^k(M,{\cal
D})$ is similar, using again Theorem \ref{2.2}.

\hfill$\Box$

In order to relate the Nambu-Poisson cohomology of a Nambu-Poisson
manifold $(M,\Lambda)$ of order $n$, $n\geq 3,$ with the foliated
cohomology of $(M,{\cal D}),$ we introduce the monomorphisms of
$C^\infty(M,\R)$-modules
\begin{equation}\label{r}
\tilde{i}^k:\Omega^k({\cal D})\rightarrow
C^{k}(\Omega^{n-1}(M)/\ker\#_{n-1};C^\infty(M,\R)),\makebox[1cm]{}[\alpha]\mapsto
\tilde{i}^k([\alpha])=\psi_\alpha,
\end{equation}
where $\psi_\alpha:\Omega^{n-1}(M)/\ker\#_{n-1}\times \dots
^{(k}\dots \times \Omega^{n-1}(M)/\ker\#_{n-1}\rightarrow
C^\infty(M,\R)$ is the mapping given by
\begin{equation}\label{n4.3'}
\psi_\alpha([\alpha_1],\dots
,[\alpha_k])=\alpha(\overline{\#_{n-1}}([\alpha_1]),\dots
,\overline{\#_{n-1}}([\alpha_k])).
\end{equation}
A direct computation, using (\ref{sc}), (\ref{partial}),
(\ref{n4.2'}), (\ref{r}) and  (\ref{n4.3'}),  proves that
\[
\tilde{i}^{k+1}\circ \tilde{d}=\tilde\partial\circ \tilde{i}^k.
\]
Hence, the mappings $\tilde{i}^k$ induce  a monomorphism between
the complexes $(\Omega^*({\cal D}),\tilde{d})$ and
$(C^{*}(\Omega^{n-1}(M)/\ker\#_{n-1};C^\infty(M,\R)),\tilde\partial).$

We will denote by
\[
\widetilde{i^k}:H^k({\cal D})\rightarrow H_{NP}^k(M)
\]
the corresponding homomorphism in cohomology.

\begin{remark}\label{r4.3}{\rm
Let $(M,\Lambda)$ be a regular  Nambu-Poisson manifold of order
$n$, with $n\geq 3.$

\medskip

$(i)$ The triple  $(\displaystyle\frac{\Lambda^{n-1}(T^*M)}{\ker
\#_{n-1}}, \lcf\;,\;\rcf^{\tilde{\;}},\widetilde{\#_{n-1}})$ is a
Lie algebroid over $M$ (see Remark \ref{r3.2}) and the Lie
algebroid cohomology is just the Nambu-Poisson cohomology.

\medskip

$(ii)$ Let ${\cal F}$ be a foliation on a manifold $N$ and
$F=\bigcup_{n\in N}{\cal F}(x)\rightarrow N$ the corresponding
vector subbundle of $TN$. Then, the mapping
\[
\pi^k:\Omega^k({\cal F})=\frac{\Omega^k(N)}{\Omega^k(N,{\cal
F})}\rightarrow C^k(\Gamma(F);C^\infty(N,\R))
\]
defined by
\[
\pi^k[\alpha](X_1,\dots ,X_k)=\alpha(X_1,\dots ,X_k),
\]
for all $[\alpha]\in \Omega^k(F)$ and $X_1,\dots ,X_k\in
\Gamma(F)$ is an isomorphism  of $C^\infty(N,\R)$-modules. This
isomorphism induces an isomorphism between the foliated cohomology
of $(N,{\cal F})$ and the Lie algebroid cohomology of
$(F,[\;,\;],i)$, $i:F\rightarrow TN$ being the natural inclusion.}
\end{remark}

Using Remarks \ref{r2.2'} and \ref{r4.3}, we deduce the
following result

\begin{theorem}
Let $(M,\Lambda)$ be a regular  Nambu-Poisson manifold of order
$n$, with $n\geq 3.$ Then, the homomorphisms of
$C^\infty(M,\R)$-modules
\[
\tilde{i}^k:\Omega^k({\cal D})
\rightarrow C^k(\frac{\Omega^{n-1}(M)}{\ker
\#_{n-1}};C^\infty(M,\R))
\]
induce an isomorphism of complexes  $\tilde{i}^*:(\Omega^*({\cal
D}), \tilde{d}) \rightarrow
(C^*(\displaystyle\frac{\Omega^{n-1}(M)}{\ker
\#_{n-1}};C^\infty(M,\R)), \tilde{\partial}).$ Thus, the
Nambu-Poisson cohomology of $M$ is isomorphic to the foliated
cohomology of $(M,{\cal D})$, that is,
\[
H^k({\cal D})\cong H_{NP}^k(M),\makebox[.5cm]{}\mbox{ for all } k.
\]
\end{theorem}


\section{A homology associated with an oriented  Nambu-Poisson manifold}
\setcounter{equation}{0}
Let $M$ be an $m$-dimensional oriented
manifold and $\nu$ be a volume form on $M.$ Denote by
$\flat_{\nu}:{\cal V}^k(M)\rightarrow \Omega^{m-k}(M)$ the
isomorphism of $C^\infty(M,\R)$-modules given by
\begin{equation}\label{5.0}
\flat_{\nu}(P)=i(P)\nu,
\end{equation}
for all $P\in {\cal V}^k(M).$

Using this isomorphism and the exterior differential $d$ we can
define a homology operator $\delta_\nu$ as follows
\begin{equation}\label{5.1}
\delta_\nu=\flat_\nu^{-1}\circ d\circ \flat_{\nu}:{\cal
V}^k(M)\rightarrow {\cal V}^{k-1}(M). \end{equation} Note that
\begin{equation}\label{5.2'}
\delta_\nu(X)=\mbox{div}_\nu X,
\end{equation}
for $X\in {\frak X}(M),$ where $\mbox{div}_\nu X$ is the
divergence of the vector field $X$ with respect to $\nu$, that is,
the $C^\infty$-real valued function on $M$ which satisfies
\begin{equation}\label{5.2''}
{\cal L}_X\nu=(\mbox{div}_\nu X)\nu.
\end{equation}

The homology associated with the complex $({\cal
V}^*(M),\delta_\nu)$ is denoted by $H_*^{\nu}(M)$ and it is dual
of the de Rham cohomology of $M$, that is,
\[
H_k^{\nu}(M)\cong H_{dR}^{m-k}(M),
\]
where $H_{dR}^*(M)$ is the de Rham cohomology of $M.$
 Therefore, $H_*^\nu(M)$ does not depend of the
chosen volume form.

In order to obtain an explicit expression of the operator
$\delta_\nu$, we will prove the following lemma which will be
useful in the sequel.

\begin{lemma}\label{l4.0}
Let $M$ be an $m$-dimensional oriented manifold and $\nu$ be a
volume form on $M.$ Then, for all $P\in {\cal V}^k(M)$ and $X\in
{\frak X}(M),$ we have
\begin{equation}\label{5.2}
{\cal L}_X \flat_\nu (P)= \flat_\nu({\cal L}_X P) +
(div_{\nu}X)\flat_\nu(P).
\end{equation}
\end{lemma}
{\it Proof:} If $k=0$ or $k=1$, relation  (\ref{5.2}) follows
using (\ref{5.0}), (\ref{5.2''}) and  the properties of the Lie
derivative operator.

Proceeding  by induction on $k$, we deduce that  (\ref{5.2}) holds
for a decomposable $k$-vector. This ends the proof.
\hfill$\Box$

Now, using this result we prove the following

\begin{proposition}\label{p5.2}
Let $M$ be an $m$-dimensional oriented manifold and $\nu$ be a
volume form on $M.$ Then
\begin{equation}\label{5.4}
i(\alpha)\delta_\nu(P)=div_\nu(i(\alpha)(P))+ (-1)^ki(d\alpha)P,
\end{equation}
for all $P\in {\cal V}^k(M)$ and $\alpha\in \Omega^{k-1}(M).$
\end{proposition}
{\it Proof:} We will proceed by induction on $k$.

If $k=1,$ (\ref{5.4}) is an immediate consequence of (\ref{5.2'})
and (\ref{5.2''}) .

Next, we will assume that (\ref{5.4}) is true for $P\in {\cal
V}^{k-1}(M)$ and $\alpha\in \Omega^{k-2}(M)$ and  we will prove
that (\ref{5.4}) also holds for a decomposable $k$-vector $P$,
\[
P=X_1\wedge \dots \wedge X_k,
\]
with $X_1,\dots ,X_k\in {\frak X}(M).$ From (\ref{5.1}),
\begin{equation}\label{0}
\begin{array}{rcl}
d(\flat_\nu(P))=d(i(X_k)(\flat_\nu(X_1\wedge \dots \wedge
X_{k-1})))&=&{\cal L}_{X_k}\flat_\nu(X_1\wedge \dots \wedge
X_{k-1}) \\&& -i(X_k)\flat_\nu(\delta_\nu(X_1\wedge \dots \wedge
X_{k-1})).
\end{array}
\end{equation}
Now, using the induction hypothesis, we have
\begin{equation}\label{1}
\begin{array}{rcl}
i(\beta)( \delta_\nu(X_1\wedge \dots \wedge
X_{k-1}))&=&\displaystyle\sum_{j=1}^{k-1}(-1)^{j+k-1}div_{\nu}(\beta(X_1,\dots
,\widehat{X_j},\dots, X_{k-1})X_j)\\&& +
(-1)^{k-1}d\beta(X_1,\dots ,X_{k-1}), \end{array}
\end{equation}
for all $\beta\in \Omega^{k-2}(M).$ Thus, one deduces that
\begin{equation}\label{2}
\begin{array}{rcl}
(-1)^{k-1}\delta_\nu(X_1\wedge \dots \wedge
X_{k-1})&\kern-10pt=&\kern-15pt\displaystyle\sum_{j=1}^{k-1}(-1)^j(div_{\nu}(X_j))X_1\wedge
\dots \wedge \widehat{X_j}\wedge \dots \wedge X_{k-1}
\\&&\kern-150pt +\displaystyle\sum_{1\leq i<j \leq
k-1}(-1)^{i+j}[X_i,X_j]\wedge X_1\wedge \dots \wedge
\widehat{X_i}\wedge \dots \wedge \widehat{X_j}\wedge \dots \wedge
X_{k-1}.
\end{array}
\end{equation}

Substituting (\ref{2}) into (\ref{0}) and using Lemma \ref{l4.0},
we obtain that
\[
\begin{array}{rcl}
(-1)^kd(\flat_\nu(P))&=&\flat_\nu\left(\displaystyle\sum_{i=1}^k(-1)^i(div_\nu
X_i)X_1\wedge \dots \wedge \widehat{X_i}\wedge \dots \wedge
X_k\right.\\[5pt]&& \left. +\displaystyle\sum_{1\leq i<j\leq
k}(-1)^{i+j}[X_i,X_j]\wedge X_1\dots \wedge
\widehat{X_i}\wedge\dots\wedge\widehat{X_j}\dots \wedge
X_k\right).
\end{array}
\]
Consequently,
\begin{equation}\label{5.4'}
\begin{array}{rcl}
(-1)^k\delta_\nu(P)&=&\displaystyle\sum_{i=1}^k(-1)^i (div_\nu
(X_i))X_1\wedge \dots \wedge \widehat{X_i}\wedge \dots \wedge
X_k\\ && \kern-40pt+\displaystyle\sum_{1\leq i<j\leq
k}^k(-1)^{i+j} [X_i,X_j]\wedge X_1\wedge \dots \wedge
\widehat{X_i}\wedge\dots\wedge\widehat{X_j}\wedge \dots  \wedge
X_k.
\end{array}
\end{equation}

On the other hand, for all $\alpha\in \Omega^{k-1}(M),$ one has
\begin{equation}\label{5.5}
\begin{array}{rcl}
(-1)^{k}div_\nu(i(\alpha)(P)) +
i(d\alpha)(P)&=&\displaystyle\sum_{i=1}^k(-1)^i\alpha (X_1,\dots
,\widehat{X_i},\dots ,X_k)div_\nu X_i\\
&&\kern-80pt+\displaystyle\sum_{1\leq i<j\leq
k}(-1)^{i+j}\alpha([X_i,X_j],X_1,\dots ,\widehat{X_i},\dots
,\widehat{X_j},\dots ,X_k).
\end{array}
\end{equation}

Therefore, from (\ref{5.4'}) and (\ref{5.5}), we conclude that
(\ref{5.4}) holds for $P=X_1\wedge \dots \wedge X_k$ and for all
$\alpha\in \Omega^{k-1}(M).$ Finally, using this result, it is
easy to prove that (\ref{5.4}) holds for all $P\in {\cal V}^k(M)$
and for all $\alpha\in \Omega^{k-1}(M).$
 \hfill$\Box$

In the following, we will describe  an interesting subcomplex of
the complex $({\cal V}^*(M), \delta_\nu)$ when $M$ is a
Nambu-Poisson manifold.

Let $(M,\Lambda)$ be an $m$-dimensional  Nambu-Poisson manifold of
order $n,$ with $n\geq 3.$ For all $k\in \{1,\dots ,n\},$ we
consider the subspace of ${\cal V}^k(M)$ given by
\[
{\cal V}^k_{t}(M,\Lambda)=\{P\in {\cal V}^k(M)/i(\alpha)(P)=0,
\mbox{ for all }\alpha \in \Omega^1(M), \;\; \alpha\in \ker\#_1\}.
\]
We will assume that ${\cal V}^0_t(M,\Lambda)=C^\infty(M,\R).$

Note that if $M$ is a regular Nambu-Poisson manifold, ${\cal
V}_t^k(M,\Lambda)$ is just the space of the $k$-vectors on $M$
which are tangent to the characteristic foliation (see Remark
\ref{r2.2'}). Thus,

\begin{lemma}\label{l5.2'}
Let $M$ be a regular Nambu-Poisson manifold of order $n$, with
$n\geq 3.$ Then
\begin{equation}\label{5.10'}
{\cal V}_t^k(M,\Lambda)=\#_{n-k}(\Omega^{n-k}(M)),
\end{equation}
for all $k\in \{0,\dots ,n\}.$
\end{lemma}

\begin{remark}
{\rm If $M$ is an arbitrary Nambu-Poisson manifold of order $n$,
with $n\geq 3$, we have that
\[
\#_{n-k}(\Omega^{n-k}(M))\subseteq {\cal V}^k_t(M,\Lambda), \mbox{
for all }k\in \{0,\dots ,n\}.
\]
However, in general, (\ref{5.10'}) does not hold as shows the
following simple example. Suppose that $M$ is an oriented manifold
of dimension $m\geq 3$ and that $\nu$ is a volume form on $M$.
Suppose also that  $f$ is a $C^\infty$-real valued function on $M$
such that $f^{-1}(0)$ is a finite subset of $M$,
$f^{-1}(0)\not=\emptyset.$ Denote by $\Lambda_\nu$ the regular
Nambu-Poisson structure induced by the volume form $\nu.$ Then,
the $m$-vector $\Lambda=f\Lambda_\nu$ defines a singular
Nambu-Poisson structure of order $m$ on $M.$ Moreover, a direct
computation proves that ${\cal V}_t^k(M,\Lambda)={\cal V}^k(M)$
for all $k\in \{0,\dots ,m\}.$ On the other hand, it is clear that
if $P\in \#_{m-k}(\Omega^{m-k}(M))$ and $x\in f^{-1}(0)$ then
$P(x)=0.$ Thus,
\[
\#_{n-k}(\Omega^{n-k}(M))\not= {\cal V}_t^k(M,\Lambda)={\cal
V}^k(M),
\]
for all $k\in \{0,\dots ,m\}.$}
\end{remark}

Next, we will prove that if $M$ is an oriented Nambu-Poisson
manifold of order $n$, with $n\geq 3,$ and $\nu$  is a volume form
on $M$ then $({\cal
V}^*_t(M,\Lambda)=\displaystyle\bigoplus_{k=1,\dots ,n}{\cal
V}_t^k(M,\Lambda))$ is a subcomplex of the complex $({\cal
V}^{*}(M),\delta_\nu).$

\begin{proposition}\label{p5.1}
Let $(M, \Lambda)$ be an oriented Nambu-Poisson manifold of order
$n$, with $n\geq 3,$ and $\nu$ be a volume form on $M.$ Then
\[
\delta_\nu({\cal V}_{t}^{k}(M,\Lambda))\subseteq {\cal
V}_t^{k-1}(M,\Lambda),
\]
for all $k\in \{1,\dots ,n\}.$
\end{proposition}
{\it Proof:} Let  $\alpha$ be an $1$-form on $M$ such that
$\alpha\in \ker\#_1.$ If $P\in {\cal V}_t^k(M,\Lambda)$ then, from
(\ref{5.4}), we have
\begin{equation}\label{5.6'}
\begin{array}{rcl}
i(\alpha)\delta_\nu(P)(\alpha_1,\dots
,\alpha_{k-2})&=&div_\nu(i(\alpha\wedge \alpha_1\wedge \dots\wedge
\alpha_{k-2})(P))\\&&
+(-1)^{k}i(d(\alpha\wedge\alpha_1\wedge\dots\wedge\alpha_{k-2})(P)),
\end{array}
\end{equation}
for all $\alpha_1,\dots, \alpha_{k-2}\in \Omega^1(M).$

Since $\alpha\in \ker\#_{1}$ and $P\in {\cal V}^k_{t}(M,\Lambda),$
we obtain that
\begin{equation}\label{5.7}
\begin{array}{rcl}
i(\alpha\wedge \alpha_1\wedge \dots
\wedge\alpha_{k-2})(P)&=&i(\alpha_1\wedge \dots \wedge
\alpha_{k-2})(i(\alpha)(P))=0,\\
i(d(\alpha\wedge\alpha_1\wedge\dots \wedge \alpha_{k-2}))(P)&=& i(
\alpha_1\wedge \dots \wedge\alpha_{k-2})(i(d\alpha)(P))\\&&-
i(d(\alpha_1\wedge \dots \wedge \alpha_{k-2}))(i(\alpha)(P))\\&=&
i(\alpha_1\wedge \dots \wedge \alpha_{k-2})(i(d\alpha)(P)).
\end{array}
\end{equation}

Next, we will see that $i(d\alpha)(P)=0,$ which proves that
$\delta_\nu(P)\in {\cal V}_t^{k-1}(M,\Lambda)$ (see (\ref{5.6'})
and (\ref{5.7})).

It is clear that the $k$-vector $P$ induces two skew-symmetric
$C^\infty(M,\R)$-linear mappings
\[
\widetilde{P}:\frac{\Omega^1(M)}{\ker\#_1}\times \dots ^{(k}\dots
\times\frac{\Omega^1(M)}{\ker \#_1}\rightarrow C^\infty(M,\R),
\]
\[
\overline{P}:\#_1(\Omega^1(M))\times \dots ^{(k}\dots \times
\#_1(\Omega^1(M))\rightarrow C ^\infty(M,\R)
\]
in such a way that
\begin{equation}\label{5.13}
P(\alpha_1,\dots, \alpha_k)=\tilde{P}([\alpha_1],\dots
,[\alpha_k])=\overline{P}(\#_1(\alpha_1),\dots ,\#_1(\alpha_k)),
\end{equation}
for all $\alpha_1,\dots ,\alpha_k\in \Omega^1(M).$ Moreover, it is
easy to prove that $\overline{P}$ is a local operator, that is, if
$U$ is an open subset of $M$ and $Q_1\in \#_1(\Omega^1(M))$ is
such that $(Q_1)_{|U}\equiv 0$ then
\[
\overline{P}(Q_1,Q_2, \dots ,Q_k)_{|U}\equiv 0,
\]
for all $Q_2,\dots , Q_k\in \#_1(\Omega^1(M)).$

Now, denote by $R$ the set of the regular points of $\Lambda$
\[
R=\{x\in M/\Lambda(x)\not=0\}.
\]
$R$ and its exterior, $\mbox{ Ext}(R),$ are open subsets of $M$.
Furthermore, it is obvious that
\[
\overline{P}(\#_1(\alpha_1),\dots
,\#_1(\alpha_k))_{|\mbox{Ext}(R)}\equiv 0
\]
for all $\alpha_1,\dots ,\alpha_k\in \Omega^1(M).$ Thus, from
(\ref{5.13}), we deduce that
\[
P(y)=0,\makebox[1cm]{} \mbox{for all } y\in \mbox{Ext}(R).
\]
This implies that
\begin{equation}\label{5.14}
i(d\alpha)(P)_{|\mbox{ Ext}(R)}\equiv 0.
\end{equation}
On the other hand, the $n$-vector $\Lambda$ induces a regular
Nambu-Poisson structure of order $n$ on $R.$ Therefore, from Lemma
\ref{l5.2'}, we obtain that there exists an $(n-k)$-form $\beta$ on
$R$ such that
\[
\#_{n-k}(\beta(y))=P(y), \makebox[1cm]{} \mbox{for all }y\in R.
\]
Consequently, if $y\in R$
\[
i(d\alpha(y))(P(y))=i(\beta(y))(\#_2(d\alpha(y))),
\]
and by  Proposition \ref{p4.0}, it follows that
\begin{equation}\label{5.14'}
i(d\alpha)(P)_{|R}\equiv 0.
\end{equation}
Finally, from (\ref{5.14}), (\ref{5.14'}) and by continuity, we
conclude that $i(d\alpha)(P)=0.$ \hfill$\Box$

Let $(M, \Lambda)$ be an  oriented Nambu-Poisson manifold of order
$n$, with $n\geq 3,$ and $\nu$ be a volume form on $M.$ Then,
Proposition \ref{p5.1} allows us to introduce the homology complex
\[
\dots \longrightarrow {\cal
V}_t^{k+1}(M,\Lambda)\stackrel{\delta_{\nu}}{\longrightarrow}{\cal
V}_t^{k}(M,\Lambda)\stackrel{\delta_\nu}{\longrightarrow}{\cal
V}_t^{k-1}(M,\Lambda)\longrightarrow\dots
\]
This complex is called the {\it canonical Nambu-Poisson complex of
$(M,\Lambda).$} The homology  of this complex is denoted by
$H^{canNP}_*(M)$ and is called {\it the canonical Nambu-Poisson
homology of} $M.$

\begin{proposition}
Let $(M, \Lambda)$ be an oriented Nambu-Poisson manifold of order
$n$, with $n\geq 3.$ The canonical Nambu-Poisson homology does not
depend on the chosen volume form.
\end{proposition}
{\it Proof:} If $\nu$ and $\nu'$ are two volume forms on $M$ then
there exists a $C^\infty$ real-valued function  $f$ on $M$ such
that $f\not=0$ at every point and
\begin{equation}\label{5.8'}
\nu'=f\nu.
\end{equation} We can
suppose, without the loss of generality, that $f>0.$

Define the  isomorphisms of $C^\infty(M,\R)$-modules
\[
\Psi^k:{\cal V}^ k_t(M,\Lambda)\rightarrow {\cal
V}_t^k(M,\Lambda)\makebox[1cm]{} P\mapsto \frac{1}{f}P,
\]
for all $k\in \{0,\dots ,n\}.$ A direct computation, using
(\ref{5.0}), (\ref{5.1})  and (\ref{5.8'}), proves that
\begin{equation}\label{516'}
\delta_{\nu'}\circ \Psi^k=\Psi^{k-1}\circ \delta_\nu.
\end{equation}
Hence, the mappings $\Psi^k$ induce an isomorphism of complexes
\[
\Psi^*:({\cal V}^*_t(M,\Lambda), \delta_{\nu})\rightarrow ( {\cal
V}_t^*(M,\Lambda),\delta_{\nu'}).
\]

\vspace{-25pt}

\hfill$\Box$

\section{Duality and the modular class of a Nambu-Poisson manifold}
\setcounter{equation}{0}

\subsection{The modular class of a Nambu-Poisson manifold}

Next, we will study when there exists a duality between the
canonical Nambu-Poisson homology and the Nambu-Poisson cohomology
of a Nambu-Poisson manifold $(M,\Lambda).$ A fundamental  tool in
this study is the modular class of $(M,\Lambda)$ which was
introduced in \cite{ILMP}. We recall its definition.

Let $(M,\Lambda)$ be an oriented $m$-dimensional Nambu-Poisson
manifold of order $n$, with $n\geq 3,$ and $\nu$ be a volume form
on $M.$

Consider the mapping ${\cal M}^\nu_{\Lambda}:C^\infty(M,\R)\times
\dots^{(n-1}\dots \times C^\infty(M,\R)\rightarrow C^\infty(M,\R)$
defined  by
\begin{equation}\label{dcm}
{\cal M}^\nu_\Lambda(f_1,\dots ,f_{n-1})=div_\nu(X_{f_1\dots
f_{n-1}}),
\end{equation}
for all $f_1,\dots ,f_{n-1}\in C^\infty(M,\R).$ Then ${\cal
M}_\Lambda^\nu$ is a skew-symmetric $(n-1)$-linear mapping and a
derivation in each argument with respect to the usual product of
functions. Thus, ${\cal M}_\Lambda^\nu$  induces an $(n-1)$-vector
on $M$ which we also denote by ${\cal M}_\Lambda^\nu$.

Moreover, the mapping
\begin{equation}\label{n6.1'}
{\cal M}_\Lambda^\nu:\Omega^{n-1}(M)\rightarrow
C^\infty(M,\R),\makebox[1cm]{} \alpha\mapsto i(\alpha){\cal
M}_\Lambda^\nu
\end{equation}
defines a $1$-cocycle in the Leibniz cohomology complex associated
with the Leibniz algebroid
$(\Lambda^{n-1}(T^*M),\lcf\;,\;\rcf,\#_{n-1})$ and its cohomology
class ${\cal M}_\Lambda=[{\cal M}^{\nu}_\Lambda]\in
H^1(\Omega^{n-1}(M);$ $C^\infty(M,\R))$ does not depend on the
chosen volume form. This cohomology class is called the {\it
modular class } of $(M,\Lambda).$

The following result proves that the $(n-1)$-vector ${\cal
M}_\Lambda^\nu$ defines also a $1$-cocycle in the Nambu-Poisson
cohomology complex.
\begin{proposition}\label{p6.1}
Let $(M,\Lambda)$ be an oriented $m$-dimensional Nambu-Poisson
manifold of order $n$, with $n\geq 3,$ and $\nu$ be a volume form
on $M.$ Then, the mapping
\begin{equation}\label{n6.1''}
\widetilde{{\cal M}^\nu_\Lambda}:\Omega^{n-1}(M)/\ker
\#_{n-1}\rightarrow C^\infty(M,\R),\makebox[1cm]{} [\alpha]\mapsto
i(\alpha){\cal M}_\Lambda^\nu,
\end{equation}
defines a $1$-cocycle in the Nambu-Poisson cohomology complex of
$(M,\Lambda)$. Moreover, the cohomology class $\widetilde{\cal
M}_\Lambda=[\widetilde{\cal M}_\Lambda^\nu]\in H_{NP}^1(M)$ does
not depend on the chosen volume form.
\end{proposition}
{\it Proof:} Let $\alpha$ be an $(n-1)$-form on $M.$ Then, using
Proposition \ref{p5.2}, we have
\begin{equation}\label{h}
div_\nu(\#_{n-1}(\alpha))=i(\alpha)\delta_\nu(\Lambda) +
(-1)^{n-1}\#_{n}(d\alpha).
\end{equation}
Now, from (\ref{n2.3'}), (\ref{dcm}) and Proposition \ref{p5.2},
it follows that
\begin{equation}\label{2*}
{\cal M}_\Lambda^\nu=\delta_\nu(\Lambda).
\end{equation}
Thus, using (\ref{h}), (\ref{2*})  and Proposition \ref{p4.0}, we
deduce that  the mapping $\widetilde{\cal M}_\Lambda^\nu$ is
well-defined.

On the other hand, since ${\cal M}_\Lambda^\nu$ defines a
$1$-cocycle in the Leibniz cohomology complex associated with the
Leibniz algebroid $(\Lambda^{n-1}(T^*M),\lcf\;,\;\rcf,\#_{n-1})$
then
\[
i(\lcf \alpha, \beta\rcf){\cal
M}_\Lambda^\nu=\#_{n-1}(\alpha)(i(\beta){\cal
M}_\Lambda^\nu)-\#_{n-1}(\beta) (i(\alpha){\cal M}_\Lambda^\nu),
\]
for all $\alpha,\beta\in \Omega^{n-1}(M).$ Therefore, we conclude
that (see (\ref{partial})),
\[
\tilde\partial \widetilde {\cal M}_{\Lambda}^\nu([\alpha],[\beta])
=\#_{n-1}(\alpha)(i(\beta){\cal M}_\Lambda^\nu)-\#_{n-1}(\beta)
(i(\alpha){\cal M}_\Lambda^\nu)-i(\lcf \alpha, \beta \rcf){\cal
M}_\Lambda^\nu=0.
\]

Finally, since the modular class of $M$ does not depend on the
chosen volume form, we deduce that the same is true for the
cohomology class $\widetilde{\cal M}_\Lambda\in H_{NP}^1(M).$
 \hfill$\Box$

\begin{remark}{\rm
Let $(M,\Lambda)$ be an oriented Nambu-Poisson manifold of order
$n$, with $n\geq 3$ and let $p^*:H^*_{NP}(M)\rightarrow
H^*(\Omega^{n-1}(M);C^\infty(M,\R))$ be the induced homomorphism
between the Nambu-Poisson cohomology of $M$ and the Leibniz
algebroid cohomology of $(\Lambda^{n-1}(T^*M),\lcf\;,\;\rcf,
\#_{n-1})$ (see Remark \ref{1-2}). Then, a direct computation,
using (\ref{n6.1'}) and (\ref{n6.1''}), proves that
\[
p^1(\widetilde{\cal M}_\Lambda)={\cal M}_\Lambda.
\]
Thus, since $p^1:H_{NP}^1(M)\rightarrow
H^1(\Omega^{n-1}(M);C^\infty(M,\R))$ is a monomorphism, it follows
that the modular class of $(M,\Lambda)$ is null if and only if
$\widetilde{\cal M}_\Lambda=0.$ }
\end{remark}

For a regular Nambu-Poisson manifold, we have

\begin{theorem}\label{63}
Let $(M,\Lambda)$ be an oriented $m$-dimensional regular
Nambu-Poisson manifold of order $n,$ with $n\geq 3.$ Then the
modular class of $(M,\Lambda)$ is null if and only if there exists
a basic volume with respect to the characteristic foliation ${\cal
D}$, that is, there exists $\mu\in \Omega^{m-n}(M)$  such that
$\mu\not=0$ at every point of $M$ and
\[
i(X_{f_{1}\dots f_{n-1}})\mu=0,\makebox[1cm]{} {\cal
L}_{X_{f_1\dots f_{n-1}}}\mu=0,
\]
for all $f_1,\dots, f_{n-1}\in C^\infty(M,\R).$
\end{theorem}
{\it Proof:} Let $\nu$ be a volume form on $M$ and suppose that
the modular class of $M$ is null. Then,  there exists $f\in
C^\infty(M,\R)$ such that
\[
{\cal M}_\Lambda^\nu=(-1)^{n-1}\#_1(df).
\]
Therefore,
\begin{equation}\label{nn1}
{\cal M}_\Lambda^\nu(df_1,\dots ,df_{n-1})=X_{f_1\dots
f_{n-1}}(f).
\end{equation}
Taking the volume form $\nu'=e^{-f}\nu$ and using (\ref{5.2''}),
(\ref{dcm}) and (\ref{nn1}), we deduce that
\begin{equation}\label{nn3}
{\cal M}_\Lambda^{\nu'}=0.
\end{equation}
Now, we consider the $(m-n)$-form
$\mu=i(\Lambda)(\nu')=\flat_{\nu'}(\Lambda)$. Then, $\mu\not=0$ at
every point of $M$ and
\[
i(X_{f_1\dots f_{n-1}})\mu=\flat_{\nu'}(\Lambda\wedge X_{f_1\dots
f_{n-1}})=0.
\]
Moreover, from (\ref{nn}), (\ref{dcm}) (\ref{nn3}) and  Lemma
\ref{l4.0} we conclude that
\[
\begin{array}{rcl}
{\cal L}_{X_{f_1\dots f_{n-1}}}\mu&=&{\cal L}_{X_{f_1\dots
f_{n-1}}}\flat_{\nu'}(\Lambda)\\ &=&\flat_{\nu'}({\cal
L}_{X_{f_1\dots f_{n-1}}}\Lambda) + (div_{\nu'} X_{f_1\dots
f_{n-1}})\flat_{\nu'}(\Lambda)=0.
\end{array}
\]

Conversely, suppose that there exists a basic volume $\mu$ with respect to
${\cal D}$. Then,
\begin{equation}\label{nn4}
i(X_{f_1\dots f_{n-1}})\mu=0,\makebox[1cm]{} {\cal L}_{X_{f_1\dots
f_{n-1}}}\mu=0,
\end{equation}
for all $f_1,\dots ,f_{n-1}\in C^\infty(M,\R).$

Let $D=\displaystyle{\cup_{x\in M}}{\cal D}(x)\rightarrow M$ be
the vector subbundle of $TM\rightarrow M$ associated with ${\cal
D}$ and $\tilde{\alpha}$ the section of the vector bundle
$\Lambda^nD^*\rightarrow M$ defined as follows. If $X_1,\dots
,X_n\in \Gamma(D)$, $\tilde\alpha(X_1,\dots ,X_n)$ is the
$C^\infty$-real valued function  on $M$ characterized by
\[
X_1\wedge \dots \wedge X_n=\tilde{\alpha}(X_1,\dots, X_n)\Lambda.
\]
Now, we extend $\tilde\alpha$ to an $n$-form $\alpha$ on $M$ such
that
\[
\alpha(X_1,\dots ,X_n)=\tilde\alpha(X_1,\dots ,X_n),
\]
for $X_1,\dots ,X_n\in \Gamma({\cal D}).$ It is clear that
\begin{equation}\label{+}
i(\Lambda)(\alpha)=1.
\end{equation}

Next, we consider the volume form $\nu$ on $M$ given by
\[
\nu=\alpha\wedge \mu.
\]
>From  (\ref{nn4}) and (\ref{+}) we have that
\begin{equation}\label{nnprima}
\flat_\nu(\Lambda)=\mu.
\end{equation}
Thus, using  (\ref{dcm}), (\ref{nn4}), (\ref{nnprima}), Lemma
\ref{l4.0} and the fact that $\mu\not=0$ at every point, we
conclude that
\[
{\cal M}_\Lambda^\nu=0.
\]
\hfill$\Box$

\begin{example}{\rm
$(i)$ Suppose that $N$ and $P$ are oriented manifolds and that
$\nu$ is a volume form on $N$. Denote by $\Lambda_{\nu}$ the
Nambu-Poisson structure on $N$ induced by the volume form $\nu$
(see Example \ref{volume}). $\Lambda_{\nu}$ defines a regular
Nambu-Poisson structure on the product manifold $M=N\times P$ and,
from Theorem \ref{63}, it follows that the modular class of
$(M,\Lambda_{\nu})$ is null.

In the same way, for a function $f\in C^{\infty}(P,\R)$ with zeros,
$f\Lambda_{\nu}$ defines a singular Nambu-Poisson structure on the
product manifold $M$ and the modular class is also null.

\medskip

$(ii)$ Let $(\frak{g}, [\, ,\, ])$ be the simple Lie algebra of dimension $3$
with basis $\{\xi,\eta,\sigma\}$ satisfying
$$[\xi,\eta]=-2\eta,\quad [\xi,\sigma]=2\sigma,\quad [\eta,\sigma]=\xi.$$

We consider a connected, simply connected, non-compact, simple Lie group $G$ such
that the Lie algebra of $G$ is $(\frak{g}, [\, ,\, ])$. From the basis
$\{\xi,\eta,\sigma\}$ one can obtain a basis of left invariant vector fields
$\{\tilde{X},\tilde{Y},\tilde{Z}\}$ on $G$ and if
$\{\tilde{\alpha},\tilde{\beta},\tilde{\gamma}\}$ is the dual basis of $1$-forms, we have that
$$d\tilde{\alpha}=\tilde{\gamma}\wedge \tilde{\beta},\quad
d\tilde{\beta}=2\tilde{\alpha}\wedge \tilde{\beta},\quad
d\tilde{\gamma}=-2\tilde{\alpha}\wedge \tilde{\gamma}.$$

Now, suppose that $S$ is a discrete subgroup such that the space $N=S\backslash
G$ of right cosets is a compact manifold (see Section $4$ of Chapter II in \cite{AGH}).
Then, the vector fields $\{\tilde{X},\tilde{Y},\tilde{Z}\}$ (respectively, the
$1$-forms $\{\tilde{\alpha},\tilde{\beta},\tilde{\gamma}\}$) induce a global
basis $\{X,Y,Z\}$ of vector fields on $N$ (respectively, a global basis
$\{\alpha,\beta,\gamma\}$ of $1$-forms on $N$) and
$$d\alpha=\gamma\wedge \beta,\quad
d\beta=2\alpha\wedge \beta,\quad
d\gamma=-2\alpha\wedge \gamma.$$

Denote by $\Lambda$ the $3$-vector on the product manifold $M=N\times S^1$ given
by
$$\Lambda=X\wedge Z\wedge E,$$
where $E$ is the dual vector field of the length element of $S^1$. It is easy to
prove that $\Lambda$ defines a regular Nambu-Poisson structure of order $3$
on $M$.

The characteristic distribution ${\cal D}$ of $(M,\Lambda)$ is the foliation
on $M$ given by $\beta=0$. Thus, ${\cal D}$ is transversally orientable
and the {\it Godbillon-Vey class} of ${\cal D}$ is the de Rham cohomology class
$4[\alpha\wedge \gamma\wedge \beta]$ (for the definition of the Godbillon-Vey class of a transversally
orientable foliation, see \cite{To} p. 29 and 30; see also \cite{GV}). It is clear
that $[\alpha\wedge \gamma\wedge \beta]\not=0$ and therefore we conclude
that it is not possible to find a basic volume with respect to ${\cal D}$ (see \cite{To},
p. 50). Consequently, from Theorem \ref{63}, we deduce that the modular class of $(M,\Lambda)$ is
not null.
}
\end{example}

\begin{remark}{\rm
Let $M$ be an oriented manifold and ${\cal D}$ an oriented
foliation on $M$ of dimension $n\geq 3$. Suppose that
$D={\bigcup_{x\in M}}{\cal D}(x)\rightarrow M$ is the vector
subbundle of $TM\rightarrow M$ associated with ${\cal D}$ and that
$\Lambda$ is a global section of the vector bundle $\Lambda^n
D\rightarrow M$, $\Lambda\not=0$ at every point. Then, $\Lambda$
defines a regular Nambu-Poisson structure of order $n$ on $M$ and
the characteristic foliation of $(M,\Lambda)$ is just ${\cal D}$.
Since $M$ is an oriented manifold, the foliation ${\cal D}$ is
transversally orientable. Thus, if the Godbillon-Vey class of
${\cal D}$ is not null, it follows that the modular class of
$(M,\Lambda)$ is not null. }
\end{remark}

\noindent

\subsection{Duality between the Nambu-Poisson cohomology and the canonical Nambu-Poisson
homology}

If $M$ is an oriented Nambu-Poisson manifold of order $n,$ with
$n\geq 3,$ and $\nu$ is a volume form on $M$, we will prove that,
under certain conditions, one can define an interesting subcomplex
of the homology complex $({\cal V}^*(M),\delta_\nu)$. In addition,
if the modular class of $M$ vanishes, we will show that there
exists a duality between the homology of this subcomplex and the
foliated cohomology of $(M,{\cal D})$, where ${\cal D}$ is the
characteristic foliation of $M$.

\begin{theorem}\label{t6.4}
Let $(M,\Lambda)$ be an oriented Nambu-Poisson manifold of order
$n$, with $n\geq 3,$ and $\nu$ be a volume form on $M.$ Then:
\begin{enumerate}
\item
$\#_*(\Omega^*(M))=\displaystyle\bigoplus_{k=0}^n(\#_{n-k}(\Omega^{n-k}(M)))$
defines a subcomplex of the homology complex $({\cal V}^*(M),$
$\delta_\nu)$ if and only if ${\cal M}_\Lambda^\nu\in
\#_1(\Omega^1(M)).$
\item If $\#_{*}(\Omega^*(M))$ is a subcomplex of
$({\cal V}^*(M),\delta_\nu)$, then the homology of this subcomplex
does not depend on the chosen volume form.
\item If the modular class of $(M,\Lambda)$ is null then
$\#_*(\Omega^*(M))$ defines a subcomplex of the homology complex $({\cal
V}^*(M),$ $ \delta_\nu)$ and
\[
\bar{H}_k^{canNP}(M)\cong H^{n-k}({\cal D}),
\]
for all $k\in \{0,\dots ,n\},$ where $H^*({\cal D})$ is the
foliated cohomology of $(M,{\cal D})$ and $\bar{H}_*^{canNP}(M)$
denotes  the homology of the complex
$(\#_*(\Omega^*(M)),\delta_\nu)$.
\end{enumerate}
\end{theorem}
{\it Proof:} $(i)$ From  (\ref{5.4}), (\ref{h}) and (\ref{2*}), we
have that
\[
\begin{array}{lcl}
i(\alpha)\delta_\nu(\#_k(\beta))&=& div_{\nu}(\#_{n-1}(\beta\wedge
\alpha))+(-1)^{n-k}\#_n(\beta\wedge d\alpha)\\ &=& i(\alpha)(
i(\beta){\cal M}_\Lambda^\nu+ (-1)^{n-1}\#_{k+1}(d\beta)),
\end{array}
\]
for all $\alpha\in \Omega^{n-k-1}(M)$ and $\beta\in \Omega^k(M).$
Thus, \begin{equation}\label{e-1}
\delta_\nu(\#_k(\beta))=(-1)^{n-1}\#_{k+1}(d\beta) + i(\beta){\cal
M}_\Lambda^\nu.
\end{equation}
Therefore, $\delta_\nu(\#_k(\Omega^k(M))\subseteq
\#_{k+1}(\Omega^{k+1}(M))$ for all $k\in \{0,\dots ,n\}$ if and
only if ${\cal M}_\Lambda^\nu\in \#_1(\Omega^1(M)).$
\medskip

$(ii)$ Let $\nu'$ be another volume form on $M$. Then, there
exists a $C^\infty$-real valued function $f$ such that $f\not=0$
at every point and $\nu'=f\nu.$ We can suppose, without the loss
of generality, that $f>0$. Thus, we can consider the isomorphisms
\[
\Psi^k:\#_k(\Omega^k(M))\rightarrow \#_k(\Omega^k(M)),
\makebox[1cm]{}P\mapsto \frac{1}{f}P.
\]
Since  $\delta_{\nu'}\circ \Psi^k=\Psi^{k-1}\circ \delta_\nu,$ it
follows that the complexes $(\#_*(\Omega^*(M)),\delta_{\nu})$ and
$(\#_*(\Omega^*(M)),$ $\delta_{\nu'})$ are isomorphic.

$(iii)$ If the modular class of $M$ is null, there exists $f\in
C^\infty(M,\R)$ such that (see (\ref{opal}))
\begin{equation}\label{e-2}\ {\cal
M}_\Lambda^\nu=\#_{1}((-1)^{n-1}df).
\end{equation}
Consequently, from  $(i),$ one deduces that
$\#_*(\Omega^*(M))$ defines a subcomplex of $({\cal
V}^*(M),\delta_\nu)$.

On the other hand, using Proposition \ref{p4.0}, we can define the
isomorphisms of $C^\infty(M,\R)$-modules
\[
h_k:\Omega^{n-k}({\cal D})=\Omega^{n-k}(M)/\ker\#_{n-k}\rightarrow
\#_{n-k}(\Omega^{n-k}(M)), \makebox[1cm]{}
h_k([\alpha])=e^{-f}\#_{n-k}(\alpha).
 \]
>From (\ref{n6.1'}), (\ref{e-1}) and (\ref{e-2}) it follows  that
$h_k\circ \tilde{d}=(-1)^{n-1}\delta_\nu\circ h_{k+1}$, where
$\tilde{d}$ is the foliated differential of $(M,{\cal D})$. So,
the above isomorphisms   induce an isomorphism  between the
cohomology group $H^{n-k}({\cal D})$  and the homology group
$\bar{H}^{canNP}_k(M).$ \hfill$\Box$

Using Remark \ref{r4.3} and Theorems \ref{63} and \ref{t6.4}, we deduce that

\begin{corollary}
Let $(M,\Lambda)$ be an oriented regular Nambu-Poisson manifold of order $n$,
with $n\geq 3$. If there exists a basic volume with respect to the characteristic
foliation ${\cal D}$ of $(M,\Lambda)$ then
$$H^k_{NP}(M)\cong H^k({\cal D})\cong H^{canNP}_{n-k}(M),$$
for all $k\in \{ 0,\ldots ,n\}$.
\end{corollary}

\section{A singular Nambu-Poisson structure}

Consider on $\R^3$ the $3$-vector defined by
\begin{equation}\label{l6.11b}
\Lambda=(x_1^2+ x_2^2 + x_3^2)\frac{\partial}{\partial x_1}\wedge
\frac{\partial}{\partial x_2}\wedge \frac{\partial}{\partial x_3},
\end{equation}
where $(x_1,x_2,x_3)$ denote the usual coordinates on $\R^3$. The
$3$-vector $\Lambda$ defines a singular Nambu-Poisson structure of
order $3$ on  $\R^3.$ Let $\nu$ be the volume form given by
\[
\nu=dx_1\wedge dx_2 \wedge dx_3.
\]
A direct computation proves that
\[
X_{x_1x_2}=(x_1^2+ x_2^2 + x_3^2)\frac{\partial}{\partial
x_3},\makebox[.5cm]{} X_{x_1x_3}=-(x_1^2+ x_2^2 +
x_3^2)\frac{\partial}{\partial x_2},\makebox[.5cm]{}
X_{x_2x_3}=(x_1^2+ x_2^2 + x_3^2)\frac{\partial}{\partial x_1},
\]
and therefore (see (\ref{dcm}))
\[ {\cal M}_\Lambda^\nu=2x_3\frac{\partial}{\partial x_1}\wedge
\frac{\partial}{\partial x_2}-2x_2\frac{\partial}{\partial
x_1}\wedge \frac{\partial}{\partial
x_3}+2x_1\frac{\partial}{\partial x_2}\wedge
\frac{\partial}{\partial x_3}.
\]
Now, if  the modular class of $(\R^3,\Lambda)$ were null then there
exists  $f\in C^\infty(\R^3,\R)$ such that
\[
i(\alpha){\cal M}_\Lambda^\nu=\#_2\alpha(f),
\]
for all $\alpha\in \Omega^2(\R^3).$ Taking the $2$-forms
$dx_1\wedge dx_2, dx_1\wedge dx_3, dx_2\wedge dx_3$ we would deduce that
\begin{equation}\label{contr}
2x_j=(x_1^2+ x_2^2+x_3^2)\frac{\partial f}{\partial
x_j},\makebox[1cm]{}\mbox{ for all }j=1,2,3.
\end{equation}
Then,
\[
f_{|\R^3-\{(0,0,0)\}}=ln(x_1^2+x_2^2+x_3^2)+c, \quad \hbox{with} \quad  c\in \R.
\]
However, this is not possible because of $f\in C^\infty(\R^3,\R).$
Thus, the modular class of $(\R^3,\Lambda)$ is not null.

Next, we will prove that there is no duality between the
Nambu-Poisson cohomology and the canonical Nambu-Poisson homology
of $(\R^3,\Lambda)$. In fact, we will show that
$$H^1_{NP}(\R^3)\not\cong H^{canNP}_2(\R^3).$$

First, we compute $H^1_{NP}(\R^3).$ In order to do this, we will
proceed  as follows:

Since $\ker\#_2=\{0\},$ then
\[
\Omega^2(\R^3)\cong \#_2(\Omega^2(\R^3))=\{(x_1^2+x_2^2+x_3^2)X/X\in
{\frak X}(\R^3)\}.
\]
This fact implies that  one can identify the co-chains
$c^1:\#_2(\Omega^2(\R^3))\rightarrow C^\infty(\R^3,\R)$ of the
Nambu-Poisson cohomology complex with the $1$-forms on $\R^3$ using the
isomorphism :
\[
\Phi:C^1(\Omega^2(\R^3);C^\infty(\R^3,\R))\rightarrow \Omega^1(\R^3),
\makebox[1cm]{} (c^1:\Omega^2(\R^3)\rightarrow C^\infty(\R^3,\R))
\mapsto \alpha
\]
such that $\alpha(X)=c^1(\beta)$, where $\#_2(\beta)=(x_1^2 +
x_2^2 + x_3^2)X.$

Under this identification the first  Nambu-Poisson cohomology
group $H^1_{NP}(\R^3)$  is the quotient space
\begin{equation}\label{l6.12b}
\frac{\{\alpha\in
\Omega^1(\R^3)/(x_1^2+x_2^2+x_3^2)d\alpha-d(x_1^2+x_2^2+x_3^2)\wedge
\alpha=0\}}{\{(x_1^2+x_2^2+x_3^2)dg/g\in C^\infty(\R^3,\R)\}}.
\end{equation}
Now, we consider the set
\[{\cal G}=\{g\in
C^\infty(\R^3-\{(0,0,0)\},\R)/(x_1^2+x_2^2+x_3^2)\frac{\partial g
}{\partial x_i}\in C^\infty(\R^3,\R),\mbox{ for all }i\in
\{1,2,3\}\}
\]
and the linear map
\[
{\cal T} :{\cal G} \rightarrow H^1_{NP}(\R^3)
\]
defined by ${\cal T}(g)=[(x_1^2+x_2^2+x_3^3)dg].$ It is clear that
the kernel of this mapping is the space $C^\infty(\R^3,\R).$
Moreover, ${\cal T}$ is an epimorphism. In fact, if $[\alpha]\in H^1_{NP}(\R^3)$,
from (\ref{l6.12b}), we deduce that in $\R^3-\{(0,0,0)\}$
$$d\left({\alpha \over x_1^2+x_2^2+x_3^2}\right)=0.$$

But this implies that there exists $g\in C^{\infty}(\R^3-\{(0,0,0)\},\R)$ such that
$${\alpha \over x_1^2+x_2^2+x_3^2}=dg$$
and therefore
$${\cal T}(g)=[\alpha].$$

Thus,
\begin{equation}\label{l6.12bb}
{{\cal G}\over{C^\infty(\R^3,\R)}}\cong H_{NP}^1(\R^3).
\end{equation}

Next, we will prove that the quotient space
$\displaystyle{{\cal G}\over {C^{\infty}(\R^3,\R)}}$ is isomorphic to $\R$.

To do that, we will use the following lemmas (a proof of the first
lemma can been found in \cite{Na}).

\begin{lemma}\cite{Na}\label{Naka}
Let $P,Q$  be two  polynomials of degree $n,$ $(n\geq 1)$ in the
indeterminates $x_1$ and $x_2$ such that satisfy
\[
(x_1^2+x_2^2)(\frac{\partial P}{\partial x_2}-\frac{\partial
Q}{\partial x_1})=2(Px_2-Qx_1).
\]
Then there exist two polynomials $\tilde{P}$,$\tilde{Q}$ of degree
$n-2$ such that $P$ and  $Q$ are written in the following form:
\[
P=ax_1 + bx_2 + (x_1^2+x_2^2)\tilde{P},\makebox[1cm]{} Q=bx_1+
ax_2 + (x_1^2+x_2^2)\tilde{Q},
\]
where $a,b$ are real constants and $\displaystyle \frac{\partial
\tilde{P}}{\partial x_2}=\displaystyle\frac{\partial
\tilde{Q}}{\partial x_1}.$
\end{lemma}

\begin{lemma}\label{ln1}
Let $A,B$ and $C$ be three polynomials of degree $n,$ $(n\geq 1)$
in the indeterminates $x_1,x_2,x_3,$  such that satisfy
\begin{equation}\label{ln}
\left.
\begin{array}{l}
(x_1^2+x_2^2+x_3^2)(\displaystyle\frac{\partial A}{\partial
x_2}-\displaystyle\frac{\partial B}{\partial
x_1})=2(Ax_2-Bx_1),\\[8pt]
(x_1^2+x_2^2+x_3^2)(\displaystyle\frac{\partial A}{\partial
x_3}-\displaystyle\frac{\partial C}{\partial
x_1})=2(Ax_3-Cx_1),\\[8pt]
(x_1^2+x_2^2+x_3^2)(\displaystyle\frac{\partial B}{\partial
x_3}-\displaystyle\frac{\partial C}{\partial x_2})=2(Ax_3-Cx_2).
\end{array}\right\}
\end{equation}
Then there exist three polynomials $\tilde{A},\tilde{B}$ and
$\tilde{C}$ of degree $n-2$ such that $A, B$ and $C$ are written
in the following form:
\[
\left.
\begin{array}{l}
A=ax_1 + (x_1^2+x_2^2+x_3^2)\tilde{A}\\ B=ax_2+
(x_1^2+x_2^2+x_3^2)\tilde{B}
\\C=ax_3+ (x_1^2+x_2^2+x_3^2)\tilde{C}
\end{array}
\right\}
\]
where $a$ is a real constant and $\displaystyle\frac{\partial
\tilde{A}}{\partial x_2}=\displaystyle\frac{\partial
\tilde{B}}{\partial x_1}$, $\displaystyle\frac{\partial
\tilde{A}}{\partial x_3}=\displaystyle\frac{\partial
\tilde{C}}{\partial x_1}$ and $\displaystyle\frac{\partial
\tilde{B}}{\partial x_3}=\frac{\partial \tilde{C}}{\partial x_2}.$
\end{lemma}
{\it Proof:} It is sufficient to prove the result for the case
when $A$, $B$ and $C$ are homogeneous polynomials. If $n=1$ is
clear that $A=ax_1$, $B=ax_2$ and $C=ax_3.$ If $n\geq 2$ we
proceed as follows.

The polynomials $A$ and $B$ can be written as
\[
A(x_1,x_2,x_3)=\sum_{k=0}^nx_3^kA_k(x_1,x_2),
\makebox[1cm]{}B(x_1,x_2,x_3)=\sum_{k=0}^nx_3^kB_k(x_1,x_2).
\]
where $A_i(x_1,x_2)$ and $B_i(x_1,x_2)$ ($i=0,\ldots ,n$) are
homogeneous polynomials in the indeterminates $x_1,x_2$.

>From the first equality of (\ref{ln}) we deduce that
\begin{equation}\label{ec1}
(x_1^2+ x_2^2)(\frac{\partial A_i}{\partial x_2}-\frac{\partial
B_i}{\partial x_1})=2(A_ix_2-B_ix_1), \makebox[.3cm]{} i\in
\{0,1\},
\end{equation}
and for all $r\in \{2,\dots n\},$
\begin{equation}\label{ec2}
(x_1^2+ x_2^2)(\frac{\partial A_{r}}{\partial
x_2}-\frac{\partial B_{r}}{\partial x_1})+ (\frac{\partial
A_{r-2}}{\partial x_2}-\frac{\partial B_{r-2}}{\partial x_1})
=2(A_rx_2-B_rx_1).
\end{equation}
Using (\ref{ec1}) and Lemma \ref{Naka} we obtain that there exist
$\tilde{A}_0,\tilde{A}_1, \tilde{B}_0$ and $\tilde{B}_1$
polynomials in the indeterminates $x_1,x_2$ such that
\[
\begin{array}{lll}
A_i=(x_1^2+x_2^2)\tilde{A}_i,\makebox[1cm]{}&
B_i=(x_1^2+x_2^2)\tilde{B}_i,\makebox[1cm]{}&
\displaystyle\frac{\partial\tilde{A}_i}{\partial
x_2}=\displaystyle\frac{\partial\tilde{B}_i}{\partial x_1},
\end{array}
\]
for $i=0,1.$

Now, from these facts and (\ref{ec2}), we have that
\[
(x_1^2+x_2^2)(\frac{\partial (A_2-\tilde{A}_0)}{\partial
x_2}-\frac{\partial (B_2-\tilde{B}_0)}{\partial
x_1})=2x_2(A_2-\tilde{A}_0)-2x_1(B_2-\tilde{B}_0).
\]
Applying again Lemma \ref{Naka} we deduce that there exist
$\tilde{A_2}$ and $\tilde{B}_2$ polynomials in the indeterminates
$x_1$ and $x_2$  such that
\[
A_2=\tilde{A}_0 + (x_1^2+x_2^2)\tilde{A}_2,\makebox[1cm]{} B_2
=\tilde{B}_0 + (x_1^2+x_2^2)\tilde{B}_2
\]
with $\displaystyle\frac{\partial \tilde{A}_2}{\partial
x_2}=\displaystyle\frac{\partial \tilde{B}_2}{\partial x_1}.$

Proceeding in a similar way we obtain a sequence of polynomials
$\tilde{A_0},\dots,  \tilde{A_n},\tilde{B_0},\dots, \tilde{B_n}$
in the indeterminates $x_1$ and $x_2$ such that
\[
A_i=(x_1^2+x_2^2)\tilde{A}_i,\makebox[1cm]{}
B_i=(x_1^2+x_2^2)\tilde{B}_i,
\]

\[
A_r=\tilde{A}_{r-2} + (x_1^2+x_2^2)\tilde{A}_r,\makebox[1cm]{}
B_r=\tilde{B}_{r-2} + (x_1^2+x_2^2)\tilde{B}_r,
\]
for $i\in \{0,1\}$ and for $r\in \{2,\dots ,n\}.$ Thus, the polynomials $A$ and $B$
can be written as
\[
A=(x_1^2+x_2^2+x_3^2)\sum_{k=0}^nx_3^k\tilde{A}_k,
\makebox[1cm]{}B=(x_1^2+x_2^2+x_3^2)\sum_{k=0}^nx_3^k\tilde{B}_k.
\]
Using the same process we also deduce that the polynomial $C$ can be written
as
\[
C=(x_1^2+x_2^2+x_3^2)\sum_{k=0}^nx_1^k\tilde{C_k},
\]
where $\tilde{C_k}$ are polynomials in the indeterminates $x_2$
and $x_3$. \hfill$\Box$

This last lemma allows us to obtain the announced result.
\begin{proposition}\label{equivo}
The quotient space $\displaystyle{{\cal G}\over{C^\infty(\R^3,\R)}}$
is isomorphic to $\R$.
\end{proposition}
{\it Proof:} Taking $g\in{\cal G}$ we have that  the $C^\infty$
real-valued functions on $\R^3$
\[
g_1=(x_1^2+x_2^2+x_3^2)\frac{\partial g}{\partial
x_1},\makebox[.5cm]{} g_2=(x_1^2+x_2^2+x_3^2)\frac{\partial
g}{\partial x_2},\makebox[.5cm]{}
g_3=(x_1^2+x_2^2+x_3^2)\frac{\partial g}{\partial x_3},
\]
satisfy
\begin{equation}\label{l6.15b}
\left.\begin{array}{l}
(x_1^2+x_2^2+x_3^2)(\displaystyle\frac{\partial g_1}{\partial
x_2}-\displaystyle\frac{\partial g_2}{\partial
x_1})=2(x_2g_1-x_1g_2),\\[8pt]
(x_1^2+x_2^2+x_3^2)(\displaystyle\frac{\partial g_1}{\partial
x_3}-\displaystyle\frac{\partial g_3}{\partial
x_1})=2(x_3g_1-x_1g_3),\\[8pt]
(x_1^2+x_2^2+x_3^2)(\displaystyle\frac{\partial g_2}{\partial
x_3}-\displaystyle\frac{\partial g_3}{\partial
x_2})=2(x_3g_2-x_2g_3).
\end{array}\right\}
\end{equation}
Then, for arbitrary $n\geq 2$, let consider the Taylor expansions
of order $n+1$ at the origin of the functions $g_1,g_2,g_3.$ We
write these Taylor expansions as $g_1=A_n + R_{1,n}$, $g_2=B_n+
R_{2,n}$ and $g_3=C_n + R_{3,n}$  where $A_n, B_n, C_n$ are
polynomials of degree $n$ which satisfy the conditions of Lemma
\ref{ln1} and $R_{i,n}$ are the remainder terms. Denote by
$[k(x_1,x_2,x_3)]_{(0,0,0)}$ the formal Taylor expansion at the
origin of $k\in C^\infty(\R^3,\R).$ Then there exists $a\in \R$
such that
\[
\left.
\begin{array}{l}
\;[g_1(x_1,x_2,x_3)-ax_1]_{(0,0,0)}=(x_1^2+x_2^2+x_3^2)A(x_1,x_2,x_3),\\
\;[g_2(x_1,x_2,x_3)-ax_2]_{(0,0,0)}=(x_1^2+x_2^2+x_3^2)B(x_1,x_2,x_3),\\
\;[g_3(x_1,x_2,x_3)-ax_3]_{(0,0,0)}=(x^2_1+x_2^2+
x_3^2)C(x_1,x_2,x_3),
\end{array}
\right\}
\]
where $A(x_1,x_2,x_3), B(x_1,x_2,x_3)$ and $C(x_1,x_2,x_3)$ are
suitable formal power series. Using Borel's theorem we have that
there exist $\alpha,\beta,\gamma\in C^\infty(\R^3,\R)$ such that
\[\left. \begin{array}{l}
\;[\alpha(x_1,x_2,x_3)]_{(0,0,0)}=A(x_1,x_2,x_3),\\
\;[\beta(x_1,x_2,x_3)]_{(0,0,0)}=B(x_1,x_2,x_3),\\
\;[\gamma(x_1,x_2,x_3)]_{(0,0,0)}=C(x_1,x_2,x_3).
\end{array}\right\}
\]
Note that the formal Taylor expansions at the origin of the
functions
\[ \alpha_1\kern-1pt=\kern-1pt g_1-ax_1-(x_1^2+ x_2^2+ x_3^2)\alpha, \makebox[.3cm]{}
\beta_1=g_2-ax_2-(x_1^2+ x_2^2+ x_3^2)\beta,\makebox[.3cm]{}
\gamma_1=g_3-ax_3-(x_1^2+ x_2^2+ x_3^2)\gamma
\]
vanish. Therefore,
 $\displaystyle\frac{\alpha_1}{(x_1^2+x_2^2+x_3^2)},$
$\displaystyle\frac{\beta_1}{(x_1^2+x_2^2+x_3^2)}$ and
$\displaystyle\frac{\gamma_1}{(x_1^2+x_2^2+x_3^2)}$ are $C^\infty$
real-valued functions on $\R^3.$

Let us consider the $C^\infty$ real-valued functions on $\R^3$
\[
h_1=\alpha + \frac{\alpha_1}{(x_1^2+x_2^2+x_3^2)},\;\;
h_2=\beta+\frac{\beta_1}{(x_1^2+x_2^2+x_3^2)},\;\;
h_3=\gamma+\frac{\gamma_1}{(x_1^2+x_2^2+x_3^2)}.
\]
Then, using (\ref{l6.15b}) and the fact that
\[
g_i=ax_i + (x_1^2+x_2^2+x_3^2)h_i,\quad i=1,2,3,
\]
we obtain that
\begin{equation}\label{h1}
\frac{\partial h_1}{\partial x_2}-\frac{\partial h_2}{\partial
x_1}=\frac{\partial h_1}{\partial x_3}-\frac{\partial
h_3}{\partial x_1}=\frac{\partial h_2}{\partial
x_3}-\frac{\partial h_3}{\partial x_2}=0.
\end{equation}

Therefore,
\begin{equation}\label{dg}
dg\kern-3pt=\kern-4pt\left(\sum_{i=1}^3\frac{g_{i}}{x_1^2+x_2^2+x_3^2}dx_i\right)
_{\kern-5pt|\R^3-\{(0,0,0)\}}\kern-20pt=
\left(d(\frac{a}{2}ln(x_1^2+x_2^2+x_3^2)) \kern-4pt+\kern-4pt
\sum_{i=1}^3h_{i}dx_i\right)_{\kern-5pt|\R^3-\{(0,0,0)\}}\kern-25pt.
\end{equation}
On the other hand, using (\ref{h1}) we deduce that
$h_1dx_1+h_2dx_2+h_3dx_3$ is a closed $1$-form on  $\R^3$ and,
since $H^1_{dR}(\R^3)=\{0\},$ we conclude that there exists
$\psi\in C^\infty(\R^3,\R)$ such that
$h_1dx_1+h_2dx_2+h_3dx_3=d\psi.$ Substituting in (\ref{dg}) we
have that
\[
g-\frac{a}{2}ln(x_1^2+x_2^2+x_3^2)=\psi_{|\R^3-\{(0,0,0)\}}+c,\quad \hbox{with}\quad
c\in \R.
\]
Consequently
\[
[g]=[\frac{a}{2}ln(x_1^2+x_2^2+x_3^2)], \makebox[1cm]{}
\mbox{with } a\in \R.
\]
This completes the proof. \hfill$\Box$

>From (\ref{l6.12bb}) and Proposition \ref{equivo}, we deduce that

\begin{proposition}
Let $\Lambda$ be the Nambu-Poisson structure on $\R^3$ given by (\ref{l6.11b}).
Then,
$$H^1_{NP}(\R^3)\cong \R.$$
\end{proposition}

On the other hand, since $\ker\#_1=\{0\}$, it follows that
\[
{\cal V}_t^k(\R^3,\Lambda)={\cal V}^k(\R^3)
\]
for all $k$. Thus, the canonical Nambu-Poisson homology of
$(\R^3,\Lambda)$ is dual of the de Rham cohomology. In particular,
$H_2^{canNP}(\R^3)\cong H^1_{dR}(\R^3)=\{0\}.$

This implies that $H^{1}_{NP}(\R^3)\not\cong H_{2}^{canNP}(\R^3)$
and therefore the duality between the Nambu-Poisson cohomology and
the canonical Nambu-Poisson homology does not hold.

\begin{remark}{\rm

$(i)$ If $\#_r:\Omega^r(\R^3)\rightarrow {\cal V}^{3-r}(\R^3)$,
$r=1,2,3$, is the induced homomorphism by the Nambu-Poisson
structure $\Lambda$ on $\R^3$, then, it is clear that $\#_r$ is a
monomorphism. Therefore, if ${\cal D}$ is the characteristic
foliation of $(\R^3,\Lambda)$, we have that the foliated
cohomology of $(\R^3,{\cal D})$ is isomorphic to the de Rham
cohomology. In particular, $H^1_{NP}(\R^3)\not\cong H^1({\cal
D})=\{0\}$. Consequently, the Nambu-Poisson cohomology and the
foliated cohomology are not isomorphic.

\medskip

$(ii)$ A direct computation shows that ${\cal
M}_\Lambda^{\nu}\not\in \#_1(\Omega^1(\R^3))$. Thus,
$\#_*(\Omega^*(\R^3))=$\linebreak
$\displaystyle\bigoplus_{k=0,\ldots ,3}\#_k(\Omega^k(\R^3))$ is
not a subcomplex of the homology complex $({\cal
V}^*(\R^3),\delta_{\nu})$ (see Theorem \ref{t6.4}).

}
\end{remark}

\section*{Acknowledgments}

This work has been partially supported through grants DGICYT
(Spain) (Projects PB97-1257 and  PB97-1487) and Project U.P.V.
127.310-EA147/98. In addition, we like to thank several
institutions for their hospitality while work on this project was
being done:  the Department of Fundamental Mathematics from
University of  La Laguna (R. Ib\'a\~nez) and Department of
Mathematics from University  of the Basque Country  (J.C. Marrero,
E. Padr\'on). We also would like to acknowledge to Nobutada
Nakanishi and Marta Macho-Stadler for helpful discussions.


\begin{thebibliography}{999}

\bibitem{AlGu} Alekseevsky, D. and  Guha, P.: On decomposability of
Nambu-Poisson tensor. Acta Math. Univ. Commenianae {\bf 65}, 1-10
(1996)

\bibitem{AGH} Auslander, L.,  Green, L. and  Hahn, F.: {\sl Flows on
homogeneous spaces}, Annals of Math. Studies 53,  Princeton Univ.
Press, 1963

\bibitem{Br} Brylinski, J.L.: A differential complex for Poisson
manifolds.  J. Diff. Geom. {\bf 28}, 93-114 (1988)

\bibitem{ChT} Chatterjee, R. and Tahktajan, L.: Aspects of classical and quantum
Nambu mechanics.  Lett. in Math. Phys. {\bf 36}, 117-126 (1996)

\bibitem{Cu} Couvier, C.: Alg\`ebres de Leibniz: d\'efinitions,
propri\'et\'es. Ann. Scient. Ec. Norm. Sup. {\bf 27}, 1-45 (1994)

\bibitem{Cz} Czachor, M.: Lie-Nambu and Beyond.  Inter. J. Theor. Phys.
{\bf 38}, n. 1,  475-500 (1999)

\bibitem{DFST} Dito, G.,  Flato, M.,  Sternheimer,D. and Tahktajan, L.: Deformation
quantization and Nambu m/echanics. Comm. Math. Phys. {\bf 183},
1-22 (1997)

\bibitem{ELW} Evens, S.,  Lu, J-H and  Weinstein, A.: Transverse measures, the
modular class and a cohomology pairing for Lie algebroids. Quart.
J. Math. Oxford, Ser.2 {\bf 50}, 417-436 (1999)

\bibitem{K} Kacimi-Alaoui, A. El: Sur la cohomologie feuillet\'ee.
Compositio Math. {\bf 49}, 195-215 (1983)

\bibitem{Ko} Koszul, J.L.: Crochet de Schoten-Nijenhuis et
cohomologie, in Elie  Cartan et les Math. d'Aujour d' Hui.
Ast\'erisque, hors s\'eries pp. 251-271 (1985)


\bibitem{Gau} Gautheron, Ph.: Some remarks concerning Nambu
mechanics. Lett. Math. Phys. {\bf 37}, 103-116 (1996)

\bibitem{GV} Godbillon, C. and  Vey, J.: Un invariant des feuilletages de codimension
un.
C. R. Acad. Sc. Paris {\bf 273}, 92-95 (1971)

\bibitem{GM} Grabowski, J. and  Marmo, G.: Remarks on Nambu-Poisson and
Nambu-Jacobi brackets. J. Phys. A: Math. Gen. {\bf 32}, 4239-4247
(1999)

\bibitem{HMM} Hector, G., Mac{\'\i}as, E. and  Saralegi, M.: Lemme de Mosser feuillet\'e
et classification des vari\'etes de Poisson r\'eguli\`eres. Publ.
Matem\`atiques {\bf 33}, 423-430 (1989)

\bibitem{ILMM} Ib\'a\~nez, R., Le\'on, M. de, Marrero, J.C. and
Mart{\'\i}n de Diego, D.: Dynamics of generalized Poisson  and
Nambu-Poisson brackets.  J. Math. Phys. {\bf 38} (5), 2332-2344
(1997).

\bibitem{ILMP} Ib\'a\~nez, R., Le\'on, de M.,  Marrero, J.C. and  Padr\'on, E.:
Leibniz algebroid associated with a Nambu-Poisson structure.  J.
Phys. A: Math. and Gen. {\bf 32},  8129-8144 (1999)

\bibitem{Lich} Lichnerowicz, A.:  Les vari\'et\'es de Poisson et leurs
alg\'ebres de Lie associ\'ees.  J. Diff. Geom.  {\bf 12}, 253-300
(1977)

\bibitem{LX} Liu, Z.J. and Xu, P.: On quadratic Poisson structures.  Letters in
Math. Phys. {\bf 26} , 33-42 (1992)


\bibitem{L1} Loday, J.L.: {\sl Cyclic Homology}, Grund. Math.
Wissen. {\bf 301}, Springer Verlag, 1992

\bibitem{L2} Loday, J.L.: Une version non commutative des
alg\'ebres de Lie: les alg\'ebres de Leibniz.  L'Enseignement
Math. {\bf 39}, 269-293 (1993)

\bibitem{LP} Loday, J.L. and  Pirashvili, T.: Universal enveloping
algebras of Leibniz and (co)-homology. Math. Ann. {\bf 296},
139-158 (1993)

\bibitem{Ma} Mackenzie, K.: {\sl Lie Groupoids and Lie algebroids in
Differential Geometry}, London Math. Soc., Lectures Notes Ser.
Vol. {\bf 124}, Cambridge Univ. Press, 1987

\bibitem{MVV} Marmo, G., Vilasi, G. and   Vinogradov, A.M.: The local
structure of $n$-Poisson and $n$-Jacobi manifolds.   J. Geom.
Phys. {\bf 25}, 141-182 (1998)


\bibitem{Mo} Molino, P.: {\sl Riemannian foliations}, Progr. Math. {\bf 73}
Birkh\"auser, Boston, 1988

\bibitem{Na} Nakanishi, N.: Poisson cohomology of plane quadratic Poisson
structures.
Publ. Res. Inst. Math. Sci. {\bf 33}, 1 , 73-89 (1997)

\bibitem{Nak} Nakanishi, N.: On Nambu-Poisson manifolds.
 Rev. Math. Phys. {\bf 10}, 499-510 (1998)

\bibitem{N} Nambu, Y.: Generalized Hamiltonian Dynamics. Phys. Reviews D
{\bf 7}, 2405-2412 (1973)

\bibitem{Ta} Takhtajan, L.: On fundations of the generalized Nambu
mechanics. Commun. Math. Phys. {\bf 160}, 295-315 (1994)

\bibitem{To} Tondeur, Ph.: {\sl Foliations on Riemannian Manifolds},
Springer-Verlag, New York, 1988

\bibitem{V1} Vaisman, I.: Variet\'et\'es riemanniennes feuillet\'ees.  Czechosl. Math. J.
{\bf 21}, 46-75 (1971)

\bibitem{V2} Vaisman, I.: {\sl Cohomology and differential forms}, M. Dekker Inc., New
York, 1973

\bibitem{V3} Vaisman, I.: {\sl Lectures on the geometry of Poisson manifolds},
Prog. Math. 118, Birkh\"auser, Basel, 1994

\bibitem{V4} Vaisman, I.: A survey on Nambu-Poisson brackets.
 Acta Math. Univ. Commenianae, Bratislava {\bf 68}, 213-243
 (1999)


\bibitem{W0} Weinstein, A.: The modular automorphism group of a Poisson
manifold. J. Geom. Phys. {\bf 23}, 379-394
 (1997)
\bibitem{W}  Weinstein, A.: Poisson geometry.  Diff. Geom. and its Appl.
{\bf 9} , 213-238 (1998)

\bibitem{Xu} Xu, P.: Gerstenhaber algebras and BV-algebras in Poisson
geometry. Commun. Math. Phys. {\bf 200}, 545-560 (1999)

\end{thebibliography}
\end{document}